%% file: main.tex
\documentclass[journal]{AIAA}

\input{package.tex}
\input{notation.tex}

\input{color.tex}
\newboolean{includefigures}

\graphicspath{{"Figures"}}
\setboolean{includefigures}{true}

\title{Statistical Analysis of the Role of Invariant Manifolds on Robust Trajectories}

\author{Amlan Sinha \footnote{Ph.D. Candidate, Mechanical and Aerospace Engineering, Princeton University, Princeton, NJ 08540.} and Ryne Beeson \footnote{Assistant Professor, Department of Mechanical and Aerospace Engineering, Princeton University, Princeton, NJ 08540.}}
\affil{Princeton University, Princeton, NJ 08840}

\begin{document}

\maketitle


\input{sections/abstract}
\input{sections/nomenclature}
\input{sections/introduction}
\input{sections/problem_formulation}
\input{sections/dynamical_model}
\input{sections/methods}
\input{sections/results_and_discussion}
\input{sections/conclusion}
\input{sections/acknowledgement}

\bibliography{main.bbl}


\end{document}

%% file: package.tex
\usepackage[utf8]{inputenc}
\usepackage{graphicx,float,placeins,tikz,onimage,caption,subcaption,xcolor}
\usepackage{physics,amsmath,mathtools,amsfonts,mathrsfs,bm}

\let\oldopenbox\openbox
\let\openbox\relax
\usepackage{amsthm}
\let\openbox\oldopenbox

\let\oldBbbk\Bbbk
\let\Bbbk\relax
\usepackage{amssymb}
\let\Bbbk\oldBbbk

\usepackage{longtable,tabularx,booktabs,multirow}
\usepackage[version=4]{mhchem}
\usepackage{siunitx}
\usepackage{textcomp}
\usepackage{lineno}
\usepackage{ifthen}
\usepackage{enumitem}
\usepackage{url}

\PassOptionsToPackage{colorinlistoftodos,prependcaption,textsize=normalsize}{todonotes}
\usepackage{todonotes}

\usepackage{listings}
\lstdefinestyle{codestyle}{
    commentstyle=\color{codegreen},
    keywordstyle=\color{blue},
    numberstyle=\tiny\color{gray},
    stringstyle=\color{codepurple},
    basicstyle=\fontfamily{pcr}\selectfont\footnotesize,
    breakatwhitespace=false,
    breaklines=true,
    captionpos=b,
    keepspaces=true,
    numbers=left,
    numbersep=5pt,
    showspaces=false,
    showstringspaces=false,
    showtabs=false,
    tabsize=2
}

%% file: notation.tex
\newcommand{\reals}{\mathbb{R}}
 
\newcommand{\rationals}{\mathbb{Q}}
\newcommand{\integers}{\mathbb{Z}}

\newcommand{\mF}{\mathcal{F}}
\renewcommand{\P}[1][P]{\mathbb{#1}}

\newcommand{\expectation}[2]{\mathbb{E}_{#1} [ {#2} ] }



%% file: color.tex
\definecolor{White}{RGB}{255, 255, 255}
\definecolor{Silver}{RGB}{192, 192, 192}
\definecolor{Grey}{RGB}{128, 128, 128}
\definecolor{Black}{RGB}{0, 0, 0}
\definecolor{Tomato}{RGB}{255, 99, 71}
\definecolor{Red}{RGB}{255, 0, 0}
\definecolor{Maroon}{RGB}{128, 0, 0}
\definecolor{Magenta}{RGB}{172, 87, 156}
\definecolor{Yellow}{RGB}{255, 255, 0}
\definecolor{Tan}{RGB}{205, 181, 145}
\definecolor{Olive}{RGB}{128, 128, 0}
\definecolor{Lime}{RGB}{0, 255, 0}
\definecolor{Green}{RGB}{0, 128, 0}
\definecolor{Chartreuse}{RGB}{160, 252, 78}
\definecolor{Aqua}{RGB}{0, 255, 255}
\definecolor{SkyBlue}{RGB}{83, 194, 236}
\definecolor{Teal}{RGB}{0, 128, 128}
\definecolor{Blue}{RGB}{0, 0, 255}
\definecolor{DodgerBlue}{RGB}{85, 154, 248}
\definecolor{RoyalBlue}{RGB}{65, 105, 225}
\definecolor{NavyBlue}{RGB}{25, 25, 107}
\definecolor{TealBlue}{RGB}{84, 113, 171}
\definecolor{Navy}{RGB}{0, 0, 128}
\definecolor{Fuchsia}{RGB}{232, 51, 244}
\definecolor{Purple}{RGB}{128, 0, 128}
\definecolor{Orange}{RGB}{242, 132, 52}

\definecolor{colunstablethreefour}{RGB}{248, 130, 83}
\definecolor{colstablethreefour}{RGB}{136, 183, 213}
\definecolor{colunstablefivesix}{RGB}{166, 135, 186}
\definecolor{colstablefivesix}{RGB}{130, 186, 126}

\definecolor{codegreen}{RGB}{0, 153, 0}
\definecolor{codepurple}{RGB}{148, 0, 209}

\definecolor{python_deep_blue}{RGB}{27, 106, 148}
\definecolor{python_deep_orange}{RGB}{228, 116, 50}

%% file: sections/abstract.tex
\begin{abstract}
As low-thrust space missions increase in prevalence, it is becoming increasingly important to design robust trajectories against unforeseen thruster outages or missed thrust events. 
Accounting for such events is particularly important in multibody systems, such as the cislunar realm, where the dynamics are chaotic and the dynamical flow is constrained by pertinent dynamical structures. 
Yet the role of these dynamical structures in robust trajectory design is unclear.
This paper provides the first comprehensive statistical study of robust and non-robust trajectories in relation to the invariant manifolds of resonant orbits in a circular restricted three-body problem. 
For both the non-robust and robust solutions analyzed in this study, the optimal subset demonstrates a closer alignment with the invariant manifolds, while the overall feasible set frequently exhibits considerable deviations.
Robust optimal trajectories shadow the invariant manifolds as closely as the non-robust optimal trajectories, and in some cases, demonstrate closer alignment than the non-robust solutions.
By maintaining proximity to these structures, low-thrust solutions are able to efficiently utilize the manifolds to achieve optimality even under operational uncertainties. 
\end{abstract}

%% file: sections/nomenclature.tex
\section*{Nomenclature}

{\renewcommand\arraystretch{1.0}
\noindent\begin{longtable*}{@{}l @{\quad=\quad} l@{}}

$J$             & objective function \\
$f$             & natural dynamics \\
$g$             & forcing dynamics \\
$\xi$           & spacecraft state \\ 
$u$             & spacecraft control \\ 
$\omega$        & random variable \\ 
$N^\dagger$     & number of segments in the reference trajectory \\ 
$N^\omega$      & number of segments in the realization trajectory \\ 
$\tau_1$        & time where a missed thrust event begins \\
$\Delta \tau$   & duration of the missed thrust event \\
$W_S^{3:4}$     & stable invariant manifold of the 3:4 resonant orbit \\
$W_U^{3:4}$     & unstable invariant manifold of the 3:4 resonant orbit \\
$W_S^{5:6}$     & stable invariant manifold of the 5:6 resonant orbit \\
$W_U^{5:6}$     & unstable invariant manifold of the 5:6 resonant orbit \\
$\mathcal{S}$   & Poincar\'{e} surface of section \\
$\mathcal{P}(x_{nr})$                     & Forward-integrated non-robust trajectory puncture point to $\mathcal{S}$ \\
$\mathcal{P}^{-1}(x_{nr})$                & Backward-integrated non-robust trajectory puncture point to $\mathcal{S}$ \\
$\mathcal{P}(x_{r})$                      & Forward-integrated robust trajectory puncture point to $\mathcal{S}$ \\
$\mathcal{P}^{-1}(x_{r})$                 & Backward-integrated robust trajectory puncture point to $\mathcal{S}$ \\
$\hat{d}_T^\mathcal{S}$                   & shortest orthogonal distance to any invariant manifold from any trajectory puncture point \\ 
$\hat{d}_A^\mathcal{S}$                   & distance along the invariant manifold to its nearest point from any trajectory puncture point \\
$\bar{d}_T^\mathcal{S} (\mathcal{W})$     & shortest orthogonal distance to the invariant manifold $\mathcal{W}$ from $x^{\mathcal{P}} \in \{\mathcal{P}(x_{nr}), \mathcal{P}^{-1}(x_{nr})\}$ \\ 
$\bar{d}_A^\mathcal{S} (\mathcal{W})$     & distance along the invariant manifold from the separatrix to its nearest point on $\mathcal{W}$ from $x^{\mathcal{P}} \in \{\mathcal{P}(x_{nr}), \mathcal{P}^{-1}(x_{nr})\}$ \\
\end{longtable*}}

%% file: sections/introduction.tex
\section{Introduction}
\label{section: introduction}

\lettrine{L}{ow-thrust} (LT) propulsion systems are becoming increasingly popular in space missions, both large strategic science missions (e.g., Hayabusa, Dawn, Hayabusa2, Bepi-Colombo, Lucy, Psyche) and small technology demonstration missions (e.g., Mars Cube One, NEA Scout, Lunar Flashlight), due to their characteristically high specific impulses which allow them to achieve a higher payload-to-propellant ratio than their impulsive counterparts. 
However, such LT missions are often susceptible to \textit{safe mode events}, which occurs if an anomalous event (e.g., impact with space debris) forces the spacecraft to depart from its nominal operating conditions causing it to enter a protective mode during which all thruster operations are switched off - if such a safe mode event coincides with a scheduled thrust arc, it results in what is known as a missed thrust event (MTE).
Due to their characteristically long thrust arcs, MTEs are quite common for LT trajectories \cite{imken_modeling_2018}.
And, unless specifically accounted for during the preliminary mission design phase, MTEs can severely impact the mission performance, and in some cases may even result in complete mission failure depending on the remaining mission time, and the available onboard fuel. 
These effects can be especially dire in missions where success depends on certain maneuvers being performed at critical junctures along the trajectory (e.g., flybys).

\subsection{Survey Of Relevant Literature}
\label{subsection: introduction: survey of relevant literature}

While robust trajectory design (i.e., the process of designing low-thrust trajectories robust to MTEs, or more concisely, \emph{robust trajectories} from here on) has garnered interest from both industry and academia in the last decade, there exists a gap in research on understanding the global geometric and topological properties of robust trajectories in relation to the dynamical structures in multibody dynamical systems.
Existing literature in the broader area of robust trajectory design can be predominantly categorized into two schools of thought - `missed thrust design' and `missed thrust analysis'.
The missed thrust design problem which refers to the problem of designing robust trajectories differs significantly from the missed thrust analysis problem which refers to the problem of identifying locations along a nominal trajectory which are most sensitive to MTEs.

Practical approaches in designing robust trajectories are more similar to the latter.
They typically involve redesigning a nominal trajectory under various missed thrust scenarios, with lower duty cycles, or with forced coast arcs in carefully chosen points along the trajectory.
Depending on the change in key performance metrics, empirical margin allocations are made for the nominal trajectory in the propellant and/or the time of flight.
Such an approach was taken with the Dawn mission to ensure that a minimum of twenty-eight days of shutdown time could be endured at all points along the nominal trajectory without significantly compromising mission objectives \cite{rayman_coupling_2007, oh_analysis_2008}.
An almost identical margin analysis has been done for the Psyche mission as well \cite{madni_missed_2020}.
Similar techniques have been explored in various studies within the literature.
For example, Laipert and Longuski investigate the trade-off between propellant reserves and schedule delays when a nominal trajectory is made robust to MTEs \cite{laipert_automated_2015}. 
Extending this work, Laipert and Imken employ similar metrics to evaluate the effects of multiple MTEs on a nominal trajectory, using a Monte Carlo approach informed by historical missed thrust data from previous missions \cite{laipert_monte_2018}.
However, due to the decoupling between the optimization of the nominal trajectory and the quantification of the uncertainty in its realization, this strategy inadvertently risks shifting the sensitivity to a different location in the redesigned trajectory, in addition to producing solutions that are glaringly sub-optimal with respect to the nominal trajectory.

Significant progress has been made in addressing the missed thrust design problem through both deterministic and stochastic approaches.
Olympio addresses the missed thrust design problem by formulating it as a two-level stochastic optimal control problem \cite{olympio_designing_2010}.
The expected thrust function, assuming a known distribution for MTEs, has been employed by Rubinsztejn et al. to address this problem \cite{rubinsztejn_designing_2021}.
More general stochastic frameworks, such as stochastic differential dynamic programming \cite{ozaki_stochastic_2018} and belief optimal control \cite{greco_robust_2022}, have also been shown to be effective in robust trajectory design. 
However, as these stochastic approaches typically model MTEs as Gaussian disturbances, they do not accurately capture their true stochastic nature.
More recent approaches explore data-driven methods in an attempt to learn the mapping between the state of the spacecraft after an MTE has occurred and the optimal control sequence going forward using neural networks \cite{rubinsztejn_neural_2020, izzo_real-time_2021} and reinforcement learning \cite{miller_low-thrust_2019, zavoli_reinforcement_2021}. 
These methods, however, only solve a local problem (i.e., small perturbations from the nominal), and are often limited in their ability to generalize to more complex gravitational environments where the LT trajectories are more sensitive to perturbations.

Further strategies in the literature address the missed thrust design problem by formulating optimization problems with constraints on the missed thrust recovery margin (i.e., the maximum amount of time a spacecraft may be allowed to coast while still being able to reach the terminal manifold once thrusters operations are resumed) \cite{olympio_deterministic_2010}.
State-of-the-art approaches extend this concept by lifting the original optimal control problem to a higher dimensional space to solve for a \textit{reference} trajectory (the path we plan to fly) simultaneously with multiple \textit{realization} trajectories (the path we may switch to should an MTE occur) from a-priori chosen points along the reference trajectory \cite{mccarty_missed_2020, venigalla_low_2020}. 
Since the reference and the realization trajectories are solved simultaneously as a single large optimization problem, it can quickly run into computational tractability issues as the number of realizations increases. 
McCarty et al. \cite{mccarty_missed_2020} therefore choose to restrict the study to a small number of realizations, whereas Venigalla et al. \cite{venigalla_low_2020} suggest an adaptive algorithm to regulate the number of realizations such that the minimum missed thrust recovery margin across all realizations remains above a threshold. 
In both studies, however, the authors note that applying the method to multibody gravitational models \textit{may be more challenging} due to the inherently chaotic nature of the underlying dynamics. 
To the best of our knowledge, the study by McCarty et al. is the only instance where the missed thrust design problem has been studied in the context of multibody gravitational environments. 
However, because they solve for a single-point solution, their study fails to elucidate the true topological properties of the robust solution space in complex dynamical environments.

Another approach to incorporate robustness is to leverage knowledge of the natural dynamical flow within the optimization framework. 
Alizadeh and Villac adopted this approach by modifying the objective function in the optimal control problem to penalize deviations from the natural dynamical flow \cite{alizadeh_sensitivity_2013}. 
However, this penalty term applies to the integral of the deviation over the entire mission duration, which may inadvertently allow for trajectories that do not consistently align with the natural dynamical flow at all times. 
With the exception of Alizadeh and Villac, current methods for designing robust trajectories generally do not approach the problem from a dynamical systems theoretic perspective. 
Even in their work, while the trajectories are encouraged to stay close to the natural dynamical flow, the precise relationship between the flow and the resulting trajectories remains unexplored, leaving a gap in fully understanding this connection.

\subsection{Contributions Of This Work}
\label{subsection: introduction: contributions of this work} 

We provide insights into the relationship between robust LT solutions and underlying dynamical structures (DS) in multibody gravitational environments with the ultimate goal of improving the robust trajectory design process through a better understanding of the geometric and topological properties of the solution space.
To do so, we first situate the missed thrust design problem within the broader class of robust control problems, which may involve a wider class of uncertainty, randomness, or stochasticity.
We consider a minimum fuel transfer between a 3:4 resonant orbit to a 5:6 resonant orbit in the Jupiter-Europa system, a problem which was originally studied by Anderson and Lo \cite{anderson_role_2009} to understand the role of invariant manifolds in non-robust LT trajectory design. 
They were able to conclude that non-robust optimal LT trajectories, without prior knowledge of the underlying structures, naturally align with pertinent invariant manifolds in the problem.
While we follow a similar approach as Anderson and Lo, we extend their work in three key aspects.
First, to rigorously analyze the relationship between trajectories and dynamical structures, we introduce distance metrics on a Poincar\'{e} surface of section, which allow for a quantitative comparison of robust and non-robust solutions, moving beyond previous qualitative assessments.
Second, we perform a detailed statistical comparison between robust and non-robust solutions, examining how their dependence on invariant manifolds evolves under varying parameters, such as the initiation and duration of the MTEs. 
By considering a family of solutions rather than focusing on a single point solution, we elucidate properties of the solution space in relation to the underlying DS.
Finally, we differentiate between feasible and optimal solutions for both robust and non-robust cases, providing key insights into how the relationship with the invariant manifolds evolves depending on whether the solutions are merely feasible or achieve optimality, highlighting significant differences in how each class of solutions leverages the dynamical structures.

\subsection{Organization Of This Paper}
\label{subsection: introduction: organization of this paper}

The paper is organized as follows.
In \S \ref{section: robust optimal control problem formulations}, we present the general robust optimal control problem, which under certain assumptions, simplifies into the missed thrust design problem we study in this paper. 
In \S \ref{section: dynamical model}, we state the circular restricted three-body model with a brief discussion to follow on the dynamical structures it exhibits, and in \S \ref{section: analysis methods}, we describe pertinent tools from dynamical systems theory along with the metrics we use to characterize the trajectory solutions with respect to the invariant manifolds.
In \S \ref{section: results and discussion}, we first present a qualitative comparison between an example robust solution with an example non-robust solution, and then present a statistical comparison between solution families in each category to uncover differences in their global properties.
Finally, we highlight the importance of this work, discuss the limitations of the current approach, and provide an outlook on future efforts in \S \ref{section: conclusion}.

%% file: sections/problem_formulation.tex
\section{Robust Optimal Control Problem Formulation}
\label{section: robust optimal control problem formulations}

We begin this section by formulating an optimal control problem that accounts for uncertainty in terminal boundary conditions and flight path constraints, randomness in system parameters, and stochastic effects. 
We start from generality so that the reader has a wider context of the missed thrust design problem considered in this paper, and so that future efforts have guidance on where to start additional extensions of the work considered in this paper with respect to additional sources of uncertainty, randomness, and stochasticity. 
After stating the general problem definition in \S \ref{subsection: robust optimal control problem formulation: the general robust formulation}, we narrow our focus to the infinite dimensional missed thrust design problem in \S \ref{subsection: robust optimal control problem formulation: the robust (MTE) formulation}, and then, a finite dimensional restriction of the missed thrust design formulation is made in \S \ref{subsection: robust optimal control problem formulation: the restricted robust (MTE) formulation}, which allows us to perform the numerical analysis of this paper. 

\subsection{The General Robust Formulation}
\label{subsection: robust optimal control problem formulation: the general robust formulation}

Let $(\Omega, \mF, (\mF_t)_{t \geq 0}, \P)$ be a filtered probability space. 
For a given random sample $\omega \in \Omega$, consider an optimal control problem where we aim to determine an extremal control solution $u^* \in \mathcal{U}$, with $\mathcal{U}$ an admissible control set, to minimize the Bolza-type cost functional,
\begin{align}
    \label{eq: general MTD (single random variable)}
    \min_{u \in \mathcal{U}} \{ J(u; \omega) \equiv \phi(\xi_1(\omega)) + & \int_0^1 \mathcal{L}(s, \xi_s(\omega), u_s) ds \ | \ \textrm{Eqs.} \ \eqref{eq: general MTD (single random variable) dynamics equation}, \eqref{eq: general MTD (single random variable) boundary conditions}, \eqref{eq: general MTD (single random variable) path constraints} \ \textrm{are satisfied} \}.
\end{align}
We consider the finite-time problem, and without loss of generality, we normalize the time interval to be $[0, 1]$. 
In Eq. \eqref{eq: general MTD (single random variable)}, $\xi$ is the solution to a stochastic differential equation driven by the control $u \in \mathcal{U}$,
\begin{align}
    \label{eq: general MTD (single random variable) dynamics equation}
    \xi_t(\omega; u) = \xi_0 (\omega) + \int_0^t f(s, \xi_s(\omega; u), \omega) ds &+ \int_0^t g(s, \xi_s(\omega; u), u_s, \omega) d\nu(s, u_s, \omega), \quad \forall t \in [0, 1], \quad \forall \omega \in \Omega,
\end{align}
taking values in a smooth manifold $\Xi$, satisfying the probabilistic initial and terminal boundary conditions,
\begin{align}
    \label{eq: general MTD (single random variable) boundary conditions}
    \P(\xi_0(\omega) \in \Xi_0) \geq 1 - \epsilon_{\Xi_0}, 
    \quad \textrm{and} \quad
    \P(\xi_1(\omega) \in \Xi_1) \geq 1 - \epsilon_{\Xi_1}, \quad \forall \omega \in \Omega,
\end{align}
and probabilistic path constraints, 
\begin{align}
    \label{eq: general MTD (single random variable) path constraints}
    \P(\varphi_i(\xi_t(\omega)) \leq 0) \geq 1 - \epsilon_\varphi, \quad \forall t \in [0, 1], \quad \forall \omega \in \Omega, \quad \forall \varphi_i \in G.
\end{align}
The numbers $\epsilon_{\Xi_0}, \epsilon_{\Xi_1}$, and $\epsilon_\varphi$ take values in the unit interval and $G$ is collection of real-valued functions. 

The optimal control problem above contains several sources of uncertainty, randomness, and stochasticity.
For example, the initial condition $\xi_0 (\omega)$ describes aleatoric uncertainty in the spacecraft state due to navigational errors. 
The drift coefficient $f(\cdot, \omega)$ represents dynamics that are independent of the control, but may also have epistemic uncertainty (e.g., uncertainty in system parameters such as location, mass, or spherical harmonics of a central gravitational body).
The dispersion coefficient $g(\cdot, \omega)$ allows for uncertainty that is dependent on the control input, and $\nu(\cdot, \omega)$ is a stochastic forcing term that may be dependent on the control process.
In the simplest case, $\nu$ could be a Lebesgue-Stieljes measure or a Brownian motion in an It\^o integral definition. 
If it is simply a Lebesgue measure, then our dynamics are for a random differential equation (e.g., uncertainty in propulsion parameters). 
Eq. \eqref{eq: general MTD (single random variable)} therefore includes both endogenous and exogenous uncertainty, and naturally accommodates aleatoric and epistemic uncertainty in the space flight problem. 

While the optimal control problem in Eq. \eqref{eq: general MTD (single random variable)} only measures the cost for a single random sample $\omega \in \Omega$, it is natural to consider an objective function that is dependent on the complete set $\Omega$, and this can be defined in a general manner by considering the space of linear functionals $\mathcal{J}$ acting on $J$. 
By the Riesz-Markov-Kakutani theorem, we can identify any element in $\mathcal{J}$ with the action of integrating $J$ against a measure. 
Because we are interested in the case where $\Omega$ is a probability sample space, it is natural for us to restrict $\mathcal{J}$ to the case where we identify it with probability measures. 
Therefore, for the general robust problem, we ultimately consider an objective function of the form, 
\begin{align}
\psi \circ J(u; \omega) \equiv \expectation{\rationals_\psi}{ \phi(\xi_1( \omega)) + \int_0^1 \mathcal{L}(s, \xi_s(\omega), u_s) ds }, 
\end{align}
for some $\psi \in \mathcal{J}$ with identifying probability measure $\rationals_\psi$. 
The deterministic case is recovered whenever the probability measure is a Dirac distribution with support on a single sample element, $\rationals_\psi = \delta_\omega$.

Therefore the general robust optimal control problem can be stated as follows. 
Given a cost functional $\psi \in \mathcal{J}$ over the sample space $\Omega$, we seek a minimizing extremal control solution $u^* \in \mathcal{U}$ to the following problem, 
\begin{align}
    \label{eq: general MTD}
    \min_{u \in \mathcal{U}} \{ \psi \circ J(u; \omega) \equiv \expectation{\rationals_\psi}{\phi(\xi_1( \omega)) + \int_0^1 \mathcal{L}(s, \xi_s(\omega), u_s) ds} \ | \ 
    \textrm{Eqs.} \ \eqref{eq: general MTD (single random variable) dynamics equation}, \eqref{eq: general MTD (single random variable) boundary conditions}, \eqref{eq: general MTD (single random variable) path constraints} \ \textrm{are satisfied}
    \}.
\end{align}

\subsection{The Robust (MTE) Formulation}
\label{subsection: robust optimal control problem formulation: the robust (MTE) formulation}

To derive the (infinite dimensional) missed thrust design problem from the general robust formulation, we now make several simplifying assumptions.
We assume that there is no stochastic forcing in Eq. \eqref{eq: general MTD (single random variable) dynamics equation}, and hence replace $d\nu(s, u_s, \omega)$ with $ds$. 
We also assume that $f$ contains no randomness, and that randomness in $g$ occurs in a very specific way. 
In particular, we introduce random times that determine whether the dispersion coefficient (or forcing function) $g$, and consequently the control input $u$, affects the state dynamics $\xi$.
To do this, we let the random sample space be identified with the unit circle (i.e., $\Omega \simeq S^1 \simeq [0, 1]$), and introduce a collection of positive strictly increasing random times $\tau \equiv \{\tau_i (\omega) \in \reals_+ \ | \ \tau_i < \tau_{i+1}, \ \forall i \in \integers_+, \omega \in \Omega \}$.
Given a control function $u$, the dynamics of an MTE trajectory is now explicitly given by, 
\begin{align}
\label{equation: system dynamics with random times}
\xi_t(\omega; u) = \xi_0(\omega) + \int_0^t f(s, \xi_s(\omega; u)) ds &+ \int_0^{\tau_1(\omega) \wedge t} g(s, \xi_s(\omega; u), u_s) ds \nonumber \\
&+ \sum_{i \in \integers_+} \int_{\tau_{2i}(\omega)}^{\tau_{2i + 1}(\omega) \wedge t} g(s, \xi_s(\omega; u), u_s) ds, \quad \forall t \in [0, 1], \quad \forall \omega \in \Omega.
\end{align}
The uncertainty due to the random times $\tau$ specifies the initiation and duration of the MTE intervals. 
In particular, the initiation of an MTE occurs when $i \in \integers_+$ is odd and the duration of that MTE is then $\tau_{i + 1} - \tau_{i}$. 
The symbol $\wedge$ in Eq. \eqref{equation: system dynamics with random times} is the minimum operator (i.e., $a \wedge b = \min(a, b)$). 
For a given realization $\omega$, if $\tau_1(\omega) > 1$, then no MTE will occur, since we are considering the finite-time problem with $t \in [0, 1]$.

We now explicitly define our choice of admissible control sets by defining the set $\mathcal{U} \equiv PC([0, 1]; \reals^n)$ for $n \in \integers_+$ to be the piecewise continuous functions on $[0, 1]$. 
Our admissible control set will be given by $\mathcal{U}^\Omega \equiv \mathcal{U}^{S^1}$ which describes the functions from $S^1$ into $\mathcal{U}$ (equivalently $\prod_{S^1} \mathcal{U}$).   
In what follows, we make the choice that $\tau_1(0) = \tau_1(1) > 1$, and hence the sample $\omega \in \{ 0, 1 \}$ will correspond to a (deterministic) non-MTE trajectory for Eq. \eqref{equation: system dynamics with random times}.
We denote this special case, when $\omega \in \{0, 1\}$, with the $\dagger$ symbol as $u^\dagger$ and refer to the state solution $\xi^\dagger$ as the \emph{reference} solution. 
For all other cases, when $\omega \in (0, 1)$, we denote the control solution as $u^\omega$ and refer to the associated state solution $\xi^\omega$ as a \emph{realization}. 

The optimal control problem for the (infinite dimensional) missed thrust design problem can now be stated as follows. We seek to find an extremal control solution $u^* \in \mathcal{U}^{\Omega}$ to the following problem, 
\begin{align}
    \label{eq: robust MTD}
    \min_{u \in \mathcal{U}^\Omega} \{ J(u^\dagger) \equiv \phi(\xi^\dagger_1) + & \int_0^1 \mathcal{L}(s, \xi^\dagger_s, u^\dagger_s) ds 
    \ | \ \textrm{Eqs.} \ \eqref{eq: robust MTD dynamics equation}, \eqref{eq: general MTD (single random variable) boundary conditions}, \eqref{eq: general MTD (single random variable) path constraints} \ \textrm{are satisfied} \},
\end{align}
where the reference and realization dynamics are given by,
\begin{align}
    \label{eq: robust MTD dynamics equation}
    \xi^\omega_t = \xi^\omega_0 + \int_0^t f(s, \xi^\omega_s) ds &+ \int_0^{\tau_1(\omega) \wedge t} g(s, \xi^\omega_s, u^\dagger_s) ds \nonumber \\
    &+ \sum_{i \in \integers_+} \int_{\tau_{2i}(\omega)}^{\tau_{2i + 1}(\omega) \wedge t} g(s, \xi^\omega_s, u^\omega_s) ds, \quad \forall t \in [0, 1], \quad \forall \omega \in \Omega.
\end{align}

\subsection{The Restricted Robust (MTE) Formulation}
\label{subsection: robust optimal control problem formulation: the restricted robust (MTE) formulation}

\ifthenelse{\boolean{includefigures}}
{
    \begin{figure}[!htb]
        \centering
        \begin{subfigure}[b]{0.45\textwidth}
            \centering
            \includegraphics[keepaspectratio, width=\textwidth]{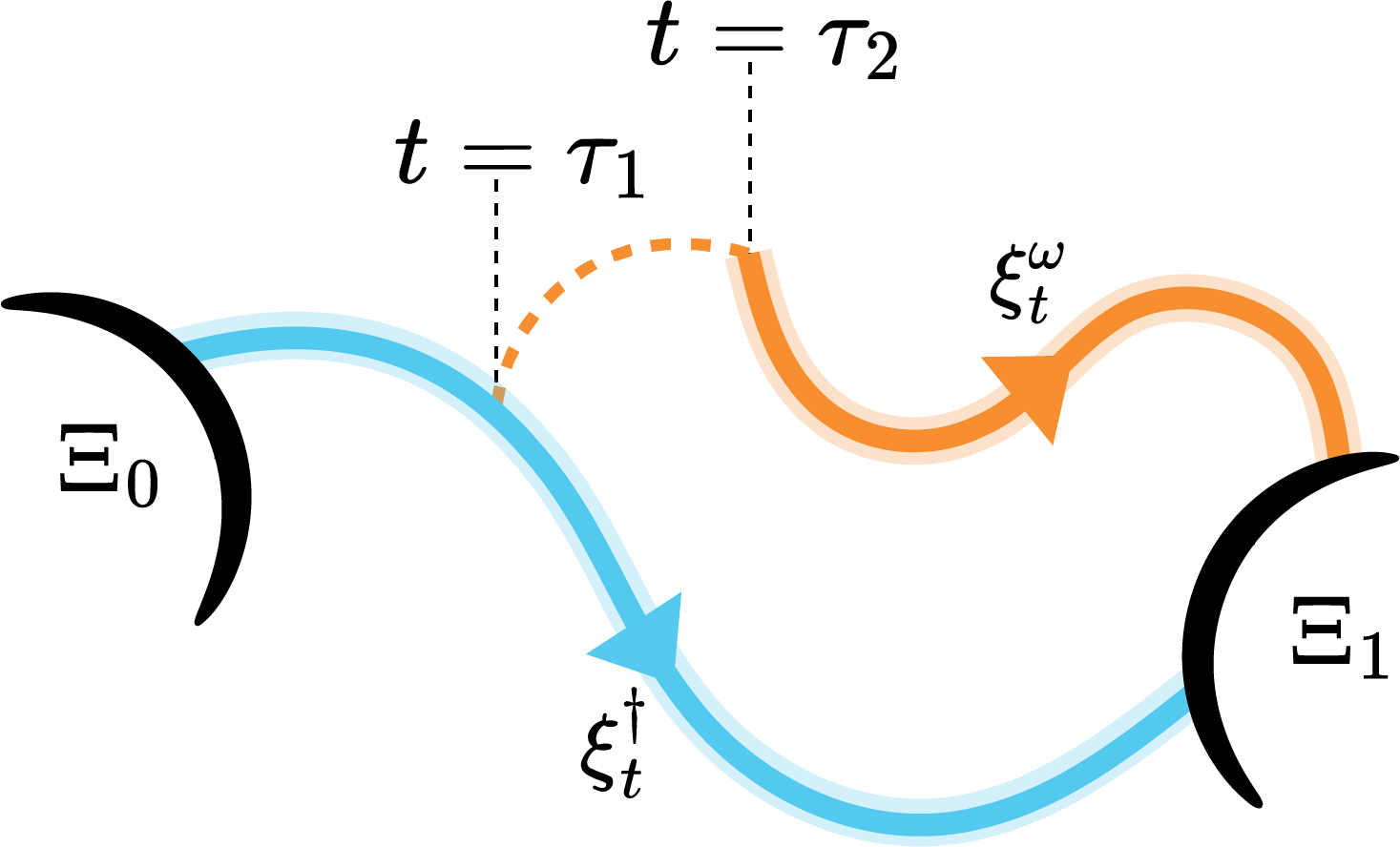}
            \caption{
            Schematic of the restricted robust problem with the reference trajectory $\xi^\dagger$ (i.e., the path we would like to fly; in \textcolor{SkyBlue}{blue}) and a realization trajectory $\xi^\omega$ (i.e., the path we may be forced to fly; in \textcolor{Orange}{orange}) is shown. 
            The initial and terminal boundary conditions for the reference solution are labeled $\Xi_0$ and $\Xi_1$ respectively.
            }
            \label{fig: problem_setup_trajectory}
        \end{subfigure}
        \hfill
        \begin{subfigure}[b]{0.45\textwidth}
            \centering
            \includegraphics[keepaspectratio, width=\textwidth]{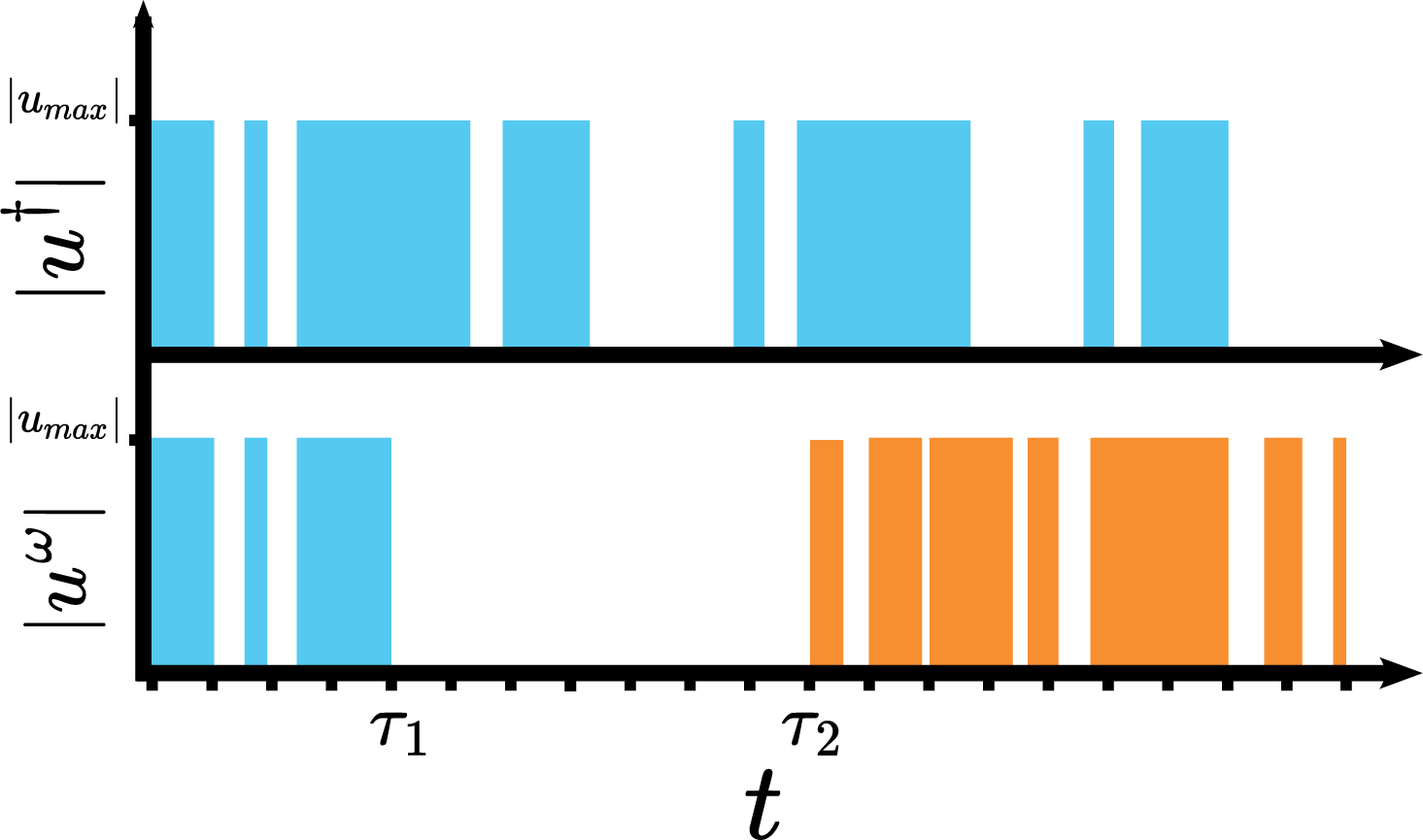}
            \caption{
            Schematic of the corresponding throttle profiles is shown.
            A thruster outage begins at $t=\tau_1$ and ends at $t=\tau_2$, which means that the spacecraft is unable to thrust during the interval $t\in[\tau_1, \tau_2]$.
            In general, the realization may exhibit a different control history compared to the reference.
            }
            \label{fig: problem_setup_throttle}
        \end{subfigure}
        \caption{
        Schematic of the restricted robust problem setup with an example robust reference trajectory along with a realization trajectory, and their corresponding throttle profiles.
        }
        \label{fig: problem_setup}
    \end{figure}
}
{
}

The main challenge in solving the robust MTE problem is the definition and approximation of the probability distribution for the random times $\tau$, which is necessary for satisfying the dynamical constraints of Eq. \eqref{eq: robust MTD dynamics equation} that couples the reference and realization solutions. 
Based on analysis of past LT missions, Imken et al. \cite{imken_modeling_2018} have suggested that the Weibull distribution is a good fit for the initiation and duration times of an MTE. 
The Weibull distribution is a continuous distribution and therefore achieving numerical tractability would require some sample approximation. 
In this paper, we assume a simpler distribution than the Weibull. 
Our assumptions are more inline with those of McCarty and Grebow \cite{mccarty_missed_2020}, and Venigalla et al. \cite{venigalla_low_2020}. 
We now describe these assumptions, by first stating them in words and then giving the mathematical definitions as subitems:

\begin{enumerate}
    \item Only one MTE will occur for any realization.
    \begin{enumerate}[label=(\roman*)]
        \item For each $\omega \in \Omega$, assume that $\tau_3(\omega) > 1$.
    \end{enumerate}
    \item Only a finite number of MTEs are allowed, with each corresponding to the start of a thrust segment (a shooting transcription is used and will be further explained in \S \ref{subsection: robust optimal control problem formulation: nonlinear program}).
    \begin{enumerate}[label=(\roman*)]
        \item Assume that $(0, 1) \subset S^1 = \Omega$ is partitioned into a collection of $N$ intervals $(E_i)_{i=1}^N$. 
        \item Assume that for every interval $E_i$, that we have for any $\omega_0, \omega_1 \in E_i$, the relation $\tau_1(\omega_0) = \tau_1(\omega_1)$.
    \end{enumerate}
    \item Only a finite number of MTE durations are allowed.
    \begin{enumerate}[label=(\roman*)]
        \item Assume that each interval $E_i$ is further partitioned into a collection of $M$ subintervals $(E_{i, j})_{j = 1}^M$. 
        \item Assume that for every subinterval $E_{i, j}$, that we have for any $\omega_0, \omega_1 \in E_{i, j}$, the relation $\tau_2(\omega_0) = \tau_2(\omega_1)$.
    \end{enumerate}
    \item Enforce deterministic boundary conditions.
    \begin{enumerate}[label=(\roman*)]
        \item Assume that $\epsilon_{\Xi_0} = \epsilon_{\Xi_1} = 0$.
    \end{enumerate}
    \item Assume that no flight path constraints exists (beyond the dynamical constraint).
    \begin{enumerate}[label=(\roman*)]
        \item Assume that $G = \varnothing$ (i.e., the empty set). 
    \end{enumerate}
    \item Lastly, we will solve $NM$ versions of our robust MTE optimal control problem, with each case corresponding to a different probability distribution $\mathbb{P}$ on $\Omega$.  
    \begin{enumerate}[label=(\roman*)]
        \item For each $i \in \{1, \hdots, N\}, j \in \{1, \hdots, M\}$, define a robust MTE optimal control problem where $\mathbb{P}(E_{i, j} \cup \{0, 1\}) = 1$.
    \end{enumerate}
\end{enumerate}

For completeness, we now state the optimal control problem for the restricted missed thrust design problem (see Fig. \ref{fig: problem_setup} for a schematic) under the additional assumptions just given:
\begin{align}
    \label{eq: simplified MTD}
    \min_{u \in \mathcal{U}^\Omega} \{ J(u^\dagger) \equiv \phi(\xi^\dagger_1) + & \int_0^1 \mathcal{L}(s, \xi^\dagger_s, u^\dagger_s) ds \ | \ \textrm{Eqs.} \ \eqref{eq: simplified MTD dynamics equation}, \ \textrm{and} \ \eqref{eq: simplified MTD boundary conditions} \ \textrm{are satisfied} \},
\end{align}
where the reference and realization dynamics are given by,
\begin{align}
    \label{eq: simplified MTD dynamics equation}
    \textrm{s.t.} \quad \xi^\omega_t = \xi^\dagger_0 + \int_0^t f(s, \xi^\omega_s) ds &+ \int_0^{\tau_1(\omega) \wedge t} g(s, \xi^\omega_s, u^\dagger_s) ds \\
    &+ \int_{\tau_2(\omega)}^{t} g(s, \xi^\omega_s, u^\omega_s) ds, \quad \forall \ t \in [0, 1], \quad \forall \omega \in \Omega, \nonumber
\end{align}
satisfying the boundary conditions,
\begin{align}
    \label{eq: simplified MTD boundary conditions}
    \xi^\omega_0 \in \Xi_0, \quad \xi^\omega_1 \in \Xi_1, \quad \forall \omega \in \Omega.
\end{align}

\subsection{Transcription to a Nonlinear Program}
\label{subsection: robust optimal control problem formulation: nonlinear program}

We solve the optimal control problem of Eq. \eqref{eq: simplified MTD} using the Dynamically Leveraged Automated (N) Multibody Trajectory Optimization (DyLAN) software developed by Beeson et al. \cite{beeson_dylan_2022}. 
DyLAN brings together dynamical systems tools with local and global optimization methods to search for solutions of optimal control problems in multibody environments. 
In this paper, we use a direct approach with a forward-backward shooting transcription to convert the optimal control problem of Eq. \eqref{eq: simplified MTD} into a nonlinear program (NLP). 
The gradient-based numerical optimizer SNOPT \cite{gill_snopt_2005} is then used to solve the NLP with initial guesses generated by the monotonic basin hopping global search algorithm \cite{wales_global_1997, leary_global_2000, englander_automated_2012, englander_tuning_2014, englander_automated_2017}.
Following the last assumption of \S \ref{subsection: robust optimal control problem formulation: the restricted robust (MTE) formulation}, we fix a version of our robust MTE problem with a probability distribution having $\mathbb{P}(E_{i, j} \cup \{0, 1\}) = 1$ for some $i \in \{1, \hdots, N\}$ and $j \in \{1, \hdots, M\}$ and then define the NLP as follows:
\begin{equation}
    \label{eq: nonlinear program}
    \begin{split}
    &\quad \min \limits_{x^\dagger \in \mathbb{R}^{\overline{N}^\dagger}, \ x^\omega \in \mathbb{R}^{\overline{N}^\omega}} \{ J(x^\dagger) = -m_f^\dagger \}, \\
    \text{subject to} &\quad c^\dagger_k(x^\dagger) = 0, \quad c^\omega_k(x^\omega) = 0, \quad \forall \ k \in \mathcal{E}, \\
    &\quad c^\dagger_k(x^\dagger) \leq 0, \quad c^\omega_k(x^\omega) \leq 0, \quad \forall \ k \in \mathcal{I},
    \end{split}
\end{equation}
where the index set $\mathcal{E}$ identifies the equality constraints, which consist of midpoint defect errors for the position, velocity, and mass continuity of the reference and realization. 
The index set $\mathcal{I}$ identifies the inequality constraints, which consists of bounds on the control variables for the reference $x^\dagger$ and realization $x^\omega$. 
The reference control decision variable has $\overline{N}^\dagger = 3 N^\dagger + 4$ components given by, 
\begin{equation}
\label{equation: decision vector}
x^\dagger = (T^\dagger_s, T^\dagger_i, T^\dagger_f, u^\dagger_1, u^\dagger_2,..., u^\dagger_{N^\dagger}, m^\dagger_f),
\end{equation}
where $T^\dagger_s$ is the shooting time, $T^\dagger_i$ the initial coast time, $T^\dagger_f$ the final coast time, and therefore the total time-of-flight is $T^\dagger_i + T^\dagger_s + T^\dagger_f$. 
$u^\dagger_k \in \mathbb{R}^3$ is a constant thrust vector for the $k$-th thrust segment that represents the throttle, in-plane, and out-of-plane thrust angle. 
Lastly, $m^\dagger_f$ is the final delivered wet mass. 
The thrust segments each have equal time of $T^\dagger_s / N^\dagger$. 

The transcription for the realization is similar, but $x^\omega$ will contain less control variables. 
The exact number is dependent on when the MTE for $x^\omega$ occurs. 
In particular, the total number of decision variables will be $\overline{N}^\omega = 3 N^\omega + 4 = 3(N^\dagger - i) + 4$, where the MTE occurs at the start of the $i$-th thrust segment for the reference solution. 
This adaptive segmentation approach, which is explained in greater detail in Sinha and Beeson \cite{amlans_lunar_2024j}, promotes congruence in control authority between the reference and realization solutions, and therefore enables a more measured understanding of the role of MTEs on the reference control solution.   
A breakdown of the number of decision variables for both the non-robust and robust cases are summarized in Table \ref{tab: decision_variables}. 

\begin{table}[hbt!]
    \caption{Number of Decision Variables (Number of Realizations = 1)}
    \label{tab: decision_variables}
    \centering
    \begin{tabular}{lcc}
        \hline
        & \textbf{Non-Robust} & \textbf{Robust} \\\hline
        \textbf{Number of Segments} & $N^\dagger$ & $N^\dagger$+$N^\omega$ \\
        \multicolumn{3}{l}{\textbf{Control Vector Components}} \\
        \hspace{3mm} \textbf{Time of Flight} & 3 & 6 \\
        \hspace{3mm} \textbf{Thrust Vector} & 3$N^\dagger$ & 3($N^\dagger$+$N^\omega$) \\
        \hspace{3mm} \textbf{Final Mass} & 1 & 2 \\
        \textbf{Number of Constraints} & 7 & 14 \\
        \hline
    \end{tabular}
\end{table}

%% file: sections/dynamical_model.tex
\section{Dynamical Model}
\label{section: dynamical model}

\subsection{Circular Restricted Three Body Problem}
\label{subsection: dynamical model: dynamical model}

In this study, we consider the motion of a spacecraft in the circular restricted three body problem (CR3BP). 
The CR3BP describes the motion of the spacecraft, whose mass is assumed to be negligible, under the influence of two celestial bodies, such as the Earth and the Moon, which rotate about their common center of mass in circular orbits.
To elucidate relevant structures in the problem, it becomes convenient to write the spacecraft's equations of motion in a synodic reference frame which rotates at the same rate as the two primaries.
The state of the spacecraft in phase space then can be described by a set of scalars $\left(q_1, q_2, q_3, \dot{q}_1, \dot{q}_2, \dot{q}_3\right)$ describing the position and the velocity.
The analysis can be further simplified by nondimensionalizing the equations using a suitable choice of units which reduces the number of parameters in the problem to one, namely the mass parameter $\mu = m_2 / (m_1 + m_2)$, where $m_1$ is the mass of the primary and $m_2 \leq m_1$ is the mass of the secondary.
With this choice of units, the gravitational constant and the mean motion both become unity and lead to the following equations of motion:
\begin{equation} 
\label{eq: eom}
\begin{split}
    \ddot{q}_1 - 2\dot{q}_2 &= -\frac{\partial}{\partial q_1} \overline{U} + \langle u, \hat{q}_1 \rangle, \\
    \ddot{q}_2 + 2\dot{q}_1 &= -\frac{\partial}{\partial q_2} \overline{U} + \langle u, \hat{q}_2 \rangle, \\
    \ddot{q}_3 &= -\frac{\partial}{\partial q_3} \overline{U} + \langle u, \hat{q}_3 \rangle,
\end{split}
\end{equation}
where $\hat{q}_i$ is the $i$-th canonical ordinate and $u$ is the control perturbation,
\begin{equation*}
    \overline{U}(r_1, r_2) \equiv -\frac{1}{2}\left((1-\mu) r_1^2 + \mu r_2^2\right) - \frac{1-\mu}{r_1} - \frac{\mu}{r_2},
\end{equation*}
is the \textit{effective} gravitational potential, 
\begin{align*}
r_1(q_1, q_2, q_3) &\equiv \sqrt{(q_1 + \mu)^2 + {q_2}^2 + {q_3}^2)}, \\
r_2(q_1, q_2, q_3) &\equiv \sqrt{(q_1 - (1-\mu))^2 + {q_2}^2 + {q_3}^2)},
\end{align*} 
are the distances between the spacecraft to the primary and the secondary respectively in the rotating frame coordinate system.

For an LT trajectory, it is also necessary to account for the change in the spacecraft mass, which can be done by simply augmenting the mass to the state of the spacecraft, where the change in the mass $m$ is governed by the differential equation:
\begin{align} 
\label{eq: mass eom}
    \dot{m} = - \frac{|u|}{I_{\textrm{sp}} g}
\end{align}
where $|u|$ is the 2-norm and hence 
the thrust magnitude, $g = 9.806 \, \text{m/s}^2$ is the gravitational acceleration on Earth and $I_{\textrm{sp}}$ is the constant specific impulse of the propulsion system. 
We neglect any other perturbations on the spacecraft (e.g., solar radiation pressure), such that the only other term perturbing the natural dynamics is the effect of the control input.

In the absence of control perturbations, there exists an integral of motion in the synodic reference frame,
\begin{equation*} 
\label{eq: jacobi integral}
    C(q_1, q_2, q_3, \dot{q}_1,\dot{q}_2,\dot{q}_3) \equiv - ({\dot{q}_1}^2 + {\dot{q}_2}^2 + {\dot{q}_3}^2) - 2\overline{U},
\end{equation*}
known as the Jacobi integral (or Jacobi constant).
The Jacobi integral, which can be thought of as a measure of the \textit{energy} of the spacecraft, remains constant between maneuvers (i.e., the Jacobi integral remains constant during the coast arcs in a spacecraft trajectory).

\subsection{Invariant Manifolds Of Unstable Periodic Orbits}
\label{subsection: dynamical model: invariant manifolds of unstable periodic orbits}

\ifthenelse{\boolean{includefigures}}
{
    \begin{figure}[!htb]
        \centering
        \begin{tikzonimage}[keepaspectratio, width=3.0 in]{"periodic_orbits.png"}
            \small
            \node[fill=white,opacity=0.85] at (0.76,0.1975){\includegraphics[keepaspectratio, width=0.75 in]{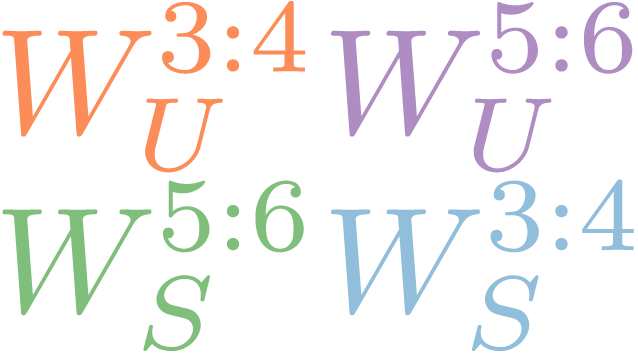}};
            \node[rectangle, draw = black, dashed, minimum width = 1.0cm, minimum height = 1.0cm, line width=0.5mm] (r1) at (.19,0.50) {};
            \node at (0.51,0.50){\includegraphics[keepaspectratio, width=0.5 in]{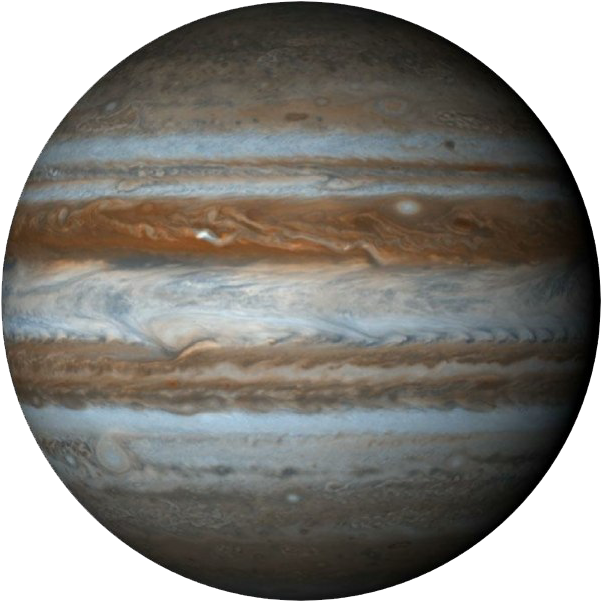}};
            \node at (1.25,0.50){\includegraphics[keepaspectratio, width=1.5 in]{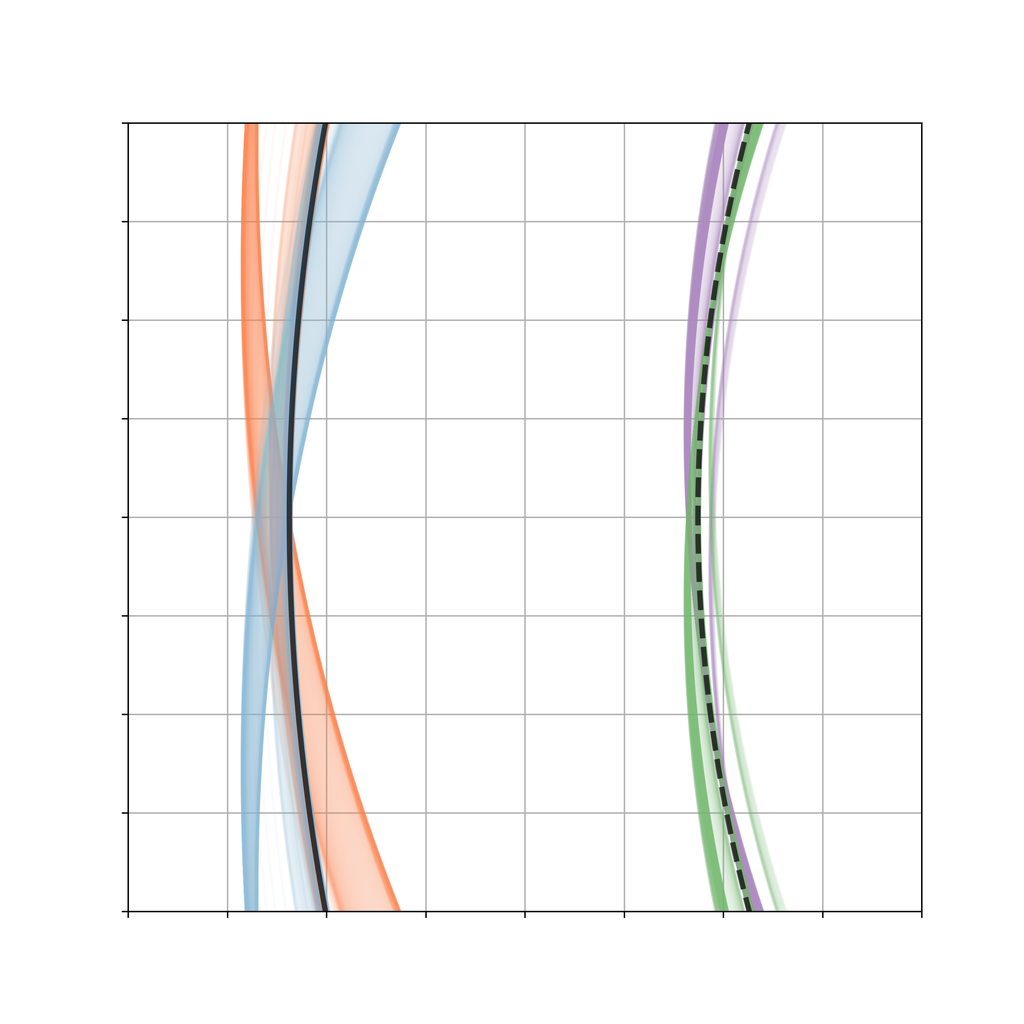}};
            \node[rectangle, draw = black, dashed, minimum width = 3.5cm, minimum height = 3.5cm, line width=0.5mm] (r2) at (1.25,0.50) {};
            \draw (r1.north west) -- (r2.north west);
            \draw[opacity=0.1] (r1.north east) -- (r2.north east);
            \draw (r1.south west) -- (r2.south west);
            \draw[opacity=0.1] (r1.south east) -- (r2.south east);
            \node[fill=white, opacity=1.0, text opacity=1, anchor=south] at (0.5,0.0) {$q_1$ [DU]};
            \node[fill=white, opacity=1.0, text opacity=1, anchor=south, rotate=90] at (0.07,0.5) {$q_2$ [DU]};
            \Large
            \node[fill=white, opacity=0.7, text opacity=1, anchor=south] at (1.135,0.45) {3:4};
            \node[fill=white, opacity=0.7, text opacity=1, anchor=south] at (1.335,0.45) {5:6};
        \end{tikzonimage}
        \caption{
        A 3:4 resonant orbit in the Jupiter-Europa CR3BP dynamical model is shown here in solid black, and a 5:6 resonant orbit is shown in dashed black. The 3:4 resonant orbit corresponds to a Jacobi integral $C_{3:4} = 2.995$, whereas the 5:6 resonant orbit correponds to a Jacobi integral $C_{3:4} = 3.005$. The 3:4 unstable invariant manifold $W_U^{3:4}$ is shown in \textcolor{colunstablethreefour}{orange}, and the corresponding stable invariant manifold $W_S^{3:4}$ is in \textcolor{colstablethreefour}{blue}. Similarly, the 5:6 unstable invariant manifold $W_U^{5:6}$ is shown in \textcolor{colunstablefivesix}{purple}, and the corresponding stable invariant manifold $W_S^{5:6}$ is in \textcolor{colstablefivesix}{green}.
        }
        \label{fig: invariant manifold}
    \end{figure}
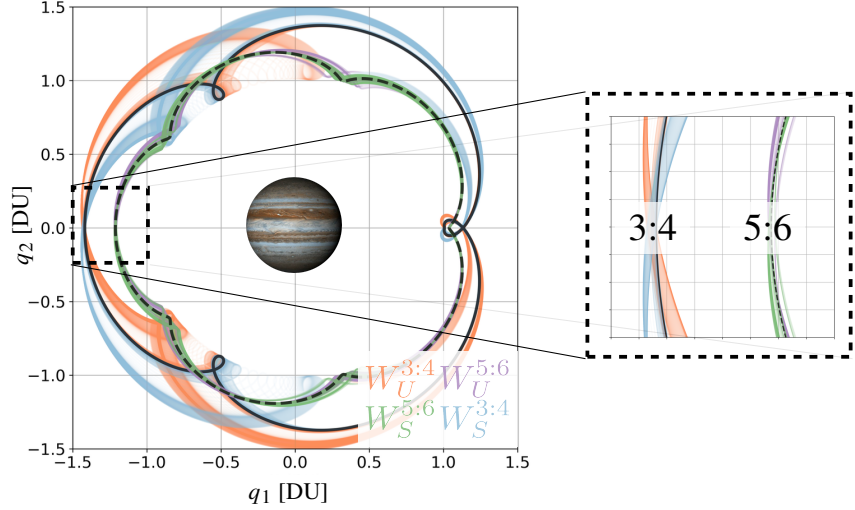
}
{
}

Despite being relatively simple, this model exhibits rich dynamical properties yielding a multitude of DS that can be leveraged during LT trajectory design. 
There exist five equilibrium points in this model known as the libration points.
Three of these points $\mathcal{L}_1, \mathcal{L}_2 \text{ and } \mathcal{L}_3$ are referred to as the colinear equilibrium points as they lie on the line joining the primary and the secondary, and the remaining two points $\mathcal{L}_4, \mathcal{L}_5$, which form equilateral triangles with the primary and the secondary as other vertices, are referred to as the non-colinear equilibrium points.
It was first shown by Poincar\'{e} \cite{poincare_methodes_1892} and later by H\'{e}non \cite{henon_exploration_1965a, henon_exploration_1965b, henon_exploration_1965c, henon_exploration_1965d, henon_exploration_1965e} that, in addition to these equilibrium points, there also exists an infinite number of periodic solutions in the three-body model.
Since then, various analytical and numerical tools have been developed to compute these periodic orbits in the neighborhood of the Lagrange points for various systems.
Each type of periodic orbit has its own distinct features, making them well-suited for specific applications.

The unstable periodic orbits within this dynamical model possess normally hyperbolic \emph{invariant manifolds}.
When we refer to a structure as ``invariant", we imply that it's time-invariant, meaning these structures remain unchanged throughout the evolution of the dynamical time.
Invariant manifolds transport material between the different realms of this model, and therefore can also be used to construct low-energy spacecraft trajectories.
They can also be classified as stable and unstable: a stable invariant manifold encompasses all points which converge to a limit set as time progresses toward infinity; conversely, an unstable manifold comprises all points which converge to the same limit set as time retrogresses toward negative infinity. 
We denote the unstable invariant manifold corresponding to a periodic orbit $\mathcal{\gamma}$ by $W^\mathcal{\gamma}_U$. 
This is the set $\mathcal{Z}$ such that $z(t) \rightarrow \mathcal{\gamma}$ as $t \rightarrow -\infty$, $\forall \ z \in \mathcal{Z}$.
Conversely, $W^\mathcal{\gamma}_S$ represents the stable invariant manifold, and comprises the set $\mathcal{Z}$ such that $z(t) \rightarrow \mathcal{\gamma}$ as $t \rightarrow \infty$, $\forall \ z \in \mathcal{Z}$.
A spacecraft originally on $\mathcal{\gamma}$ will shadow $W^\mathcal{\gamma}_S$ when perturbed in the direction of the stable eigenvector of the monodromy matrix corresponding to $\mathcal{\gamma}$, and will shadow $W^\mathcal{\gamma}_U$ when perturbed in the direction of the unstable eigenvector.
Invariant manifolds share the same Jacobi integral as the periodic orbits to which they are associated.
The invariant manifolds associated with the 3:4 and 5:6 resonant orbits in the Jupiter-Europa system are shown in Fig. \ref{fig: invariant manifold}.

The invariant manifolds describe the local dynamics in the neighborhood of the periodic orbits, and provide a global template for LT trajectories.
Previous studies have shown that optimization algorithms applied to minimum-fuel problems and without explicit prior knowledge of the underlying DS, qualitatively converge to locally optimal solutions that align themselves with these structures \cite{lo_interplanetary_2002}.
Anderson and Lo \cite{anderson_role_2009} extended the work of Lo \cite{lo_interplanetary_2002} to investigate a minimum-fuel LT moon tour in the Jupiter-Europa CR3BP, originally developed by Lam et al. \cite{lam_ganymede_2004} using the trajectory design tool Mystic \cite{whiffen_mystic_2006}. Anderson and Lo discovered that the numerical optimal trajectories indeed appear to shadow the invariant manifolds of resonant orbits. 
The main purpose of this paper is to further extend Anderson and Lo's work beyond qualitative understanding to a quantitative one, as well as to study robust trajectories and their dependence on the underlying DS.
We aim to compare the behavior of robust and non-robust trajectories in relation to these structures.
Having knowledge of the relationship of the robust solutions to the DS can be useful in developing good initial guesses and efficient algorithms for global design of robust optimal LT trajectories.

%% file: sections/methods.tex
\section{Analysis Methods}
\label{section: analysis methods}

\subsection{Poincar\'{e} Surface of Sections}
\label{subsection: analysis methods: poincare surface of sections}

\ifthenelse{\boolean{includefigures}}
{
    \begin{figure}[!htb]
    \centering
    \begin{subfigure}[b]{0.45\textwidth}
        \centering
        \begin{tikzonimage}[keepaspectratio, width=3 in]{"poincare_section.png"}
            \Large
            \node[anchor=south] at (0.58,0.58) {$x$};
            \node[anchor=south] at (0.725,0.75) {$\mathcal{S}$};
            \large
            \node[anchor=south] at (0.5,0.4) {$\mathcal{P}^{-1}(x)$};
            \node[anchor=south] at (0.7,0.55) {$\mathcal{P}(x)$};
        \end{tikzonimage}
        \caption{
        Schematic of a Poincar\'{e} map, denoted as $\mathcal{P}:\mathcal{S}\rightarrow\mathcal{S}$, obtained by intersecting a point $x$ along the trajectory propagated using the natural dynamics with a Poincar\'{e} surface of section $\mathcal{S}$.
        The forward-in-time integrated map of $x$ is represented by $\mathcal{P}(x)$ and is indicated with a circle. 
        Conversely, the backward-in-time integrated map of $x$ is represented by $\mathcal{P}^{-1}(x)$ and is indicated with a square.
        The forward-integrated map coincides with the backward-integrated map of a trajectory only if the point $x$ belongs to a periodic orbit.
        }
        \label{fig: poincare section}
    \end{subfigure}
    \hfill
    \begin{subfigure}[b]{0.45\textwidth}
        \centering
        \begin{tikzonimage}[keepaspectratio, width=2.5 in]{"invariant_manifold_example.png"}
        \node[fill=white,opacity=0.85] at (0.8275,0.14){\includegraphics[keepaspectratio, width=0.75 in]{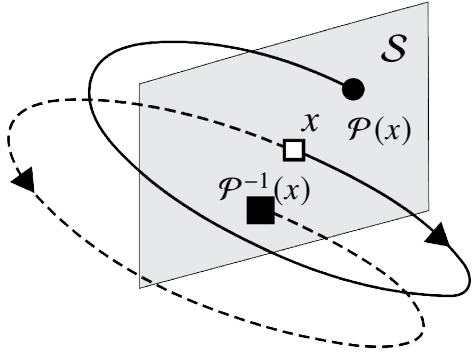}
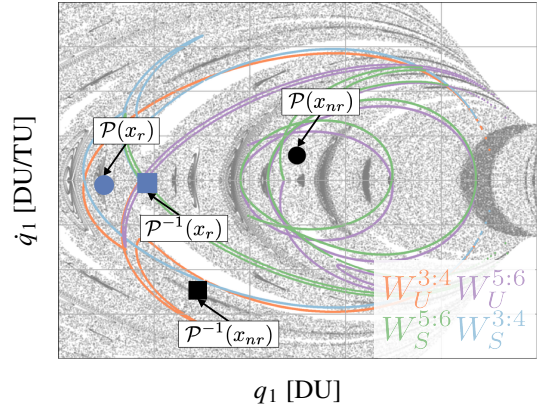};
        \node[fill=white, opacity=1.0, text opacity=1, anchor=south] at (0.5,-0.15) {$q_1$ [DU]};
        \node[fill=white, opacity=1.0, text opacity=1, anchor=south, rotate=90] at (-0.025,0.5) {$\dot{q}_1$ [DU/TU]};
        \end{tikzonimage}
        \caption{
        An example Poincar\'{e} surface of section $\mathcal{S}$, featuring invariant manifold puncture points that correspond to a 3:4 resonant orbit and a 5:6 resonant orbit in the Jupiter-Europa system as well the trajectory puncture points for an example nonrobust (\textcolor{Black}{black}) and robust (\textcolor{TealBlue}{blue}) solution.
        Also displayed are the background points (\textcolor{Grey}{grey}) that provide a visualization of the global dynamical template at a particular Jacobi integral.
        Every puncture point on a given Poincar\'{e} section possesses the same Jacobi integral.
        }
        \label{fig: poincare example}
    \end{subfigure}
    \caption{
        Schematic of a Poincar\'{e} Map along with an example Poincar\'{e} surface of section.
    }
    \label{fig: poincare section schematic}
    \end{figure}
}
{
}

Poincar\'{e} return maps, or simply Poincar\'{e} maps, are an effective tool for analyzing rotational flows such as periodic or quasi-periodic orbits, or even flow in the vicinity of a periodic orbit and can therefore be used to investigate the trajectories as well as pertinent invariant manifolds.
We consider a point $x\in\mathcal{S}$ on the surface $\mathcal{S}$ which we evolve in time according to the governing dynamical equations until it intersects $\mathcal{S}$ again transversely. 
We denote the intersection of the point $x$ with $\mathcal{S}$ as $\mathcal{P}(x)$.
Therefore, $\mathcal{P}(x)$ represents the first return of the trajectory to $\mathcal{S}$, $\mathcal{P}^2(x)$ represents the second return of the trajectory to $\mathcal{S}$ and so on. 
We can continue to evolve $x$ in time and record its state after every intersection with $\mathcal{S}$ and by doing so, we effectively reduce the global orbit structure governed by differential equations to a discrete-time dynamical system given by the map $\mathcal{P}$.
A Poincar\'{e} map can therefore be mathematically described by $\mathcal{P}$: $\mathcal{S} \rightarrow \mathcal{S}$, where $\mathcal{S}$ is referred to as the Poincar\'{e} surface of section, or simply Poincar\'{e} section (see Fig. \ref{fig: poincare section} for a visual representation).
If a trajectory is (sufficiently) planar, the map $\mathcal{P}$ provides sufficient information to fully characterize the trajectory.

An example Poincar\'{e} section is shown in Fig. \ref{fig: poincare section schematic}. 
Puncture points corresponding to the evolution of the forward integrated non-robust solution $\mathcal{P}(x_{nr})$ and robust solution $\mathcal{P}(x_{r})$ are denoted by black and blue circles respectively which were computed by mapping points along the state forward in time without thrust until it intersected $\mathcal{S}$.
Similarly, puncture points corresponding to the evolution of the backward integrated non-robust $\mathcal{P}^{-1}(x_{nr})$ and robust trajectory $\mathcal{P}^{-1}(x_{r})$ are denoted by black and blue squares.
In this study, we only record the first puncture point for trajectories, discarding the subsequent ones.

The puncture points corresponding to the stable and unstable invariant manifolds for the relevant resonant orbits are an important part of the analysis in this study.
In our study, the sets containing these puncture points are labeled as $W_i^\gamma$, where the superscript $\gamma$ allows us to discern what type of periodic orbit we are referring to and the subscript $i$ allows us to discern the stability of the invariant manifolds.
To compute the invariant manifold puncture points, one begins by considering perturbations parallel (and anti-parallel) to the eigenvectors of the monodromy matrix, which are then propagated forward (if unstable) or backward (if stable) in time until their first intersection with $\mathcal{S}$.
The states of these puncture points are recorded and plotted using a coordinate system of choice (e.g., $q_1-\dot{q}_1$). 
In this study, we consider a range of perturbation from [$1 \times 10^{-6}$, $3 \times 10^{-1}$] (in nondimensional units) are considered both parallel to and anti-parallel to the eigenvectors of the monodromy matrix, which are then propagated until their first intersection with $\mathcal{S}$.
$W_{\mathrm{U}}^{3:4}$ represents the intersection of the unstable manifold of the 3:4 resonant orbit (\textcolor{colunstablethreefour}{orange}), $W_{\mathrm{S}}^{3:4}$ represents the intersection of the stable manifold of the 3:4 resonant orbit (\textcolor{colstablethreefour}{blue}), $W_{\mathrm{U}}^{5:6}$ represents the intersection of the unstable manifold of the 5:6 resonant orbit (\textcolor{colunstablefivesix}{purple}), and $W_{\mathrm{S}}^{5:6}$ represents the intersection of the stable manifold of the 5:6 resonant orbit (\textcolor{colstablefivesix}{green}).
Approximately 10,000 puncture points are computed for each invariant manifold in the subsequent analysis.

The background points (\textcolor{Grey}{grey}) are computed by first considering points on a uniform grid on the $x$-axis on $\mathcal{S}$ and then integrating these points forward in time until they intersect the surface a number of times.
The number of points retained for subsequent analysis is largely problem dependent but should be chosen such that a sufficiently detailed visual representation of the global dynamical flow template is visible.
In this study, the background points are computed by considering $10, 000$ points on an equally spaced grid on the $x$-axis on $\mathcal{S}$ and integrating forward in time until they intersect the surface 10 times.
The first five puncture points were discarded to remove the distortion of the grid during integration, and the remaining puncture points were recorded.
The number of intersections being retained was arbitrarily chosen, but for the context of this research, this choice does not impact the results since the background points are only for visual purposes and do not affect any of the subsequent analyses.
50,000 points were recorded for the background points in total.

\subsection{Jacobi Integral of Low-Thrust Trajectories}
\label{subsection: analysis_methods: jacobi_integral_of_lt_trajectories}

\ifthenelse{\boolean{includefigures}}
{
    \begin{figure}[!htbp]
        \centering
        \begin{tikzonimage}[keepaspectratio, width=3.0 in]{"orbit_by_energy.png"}
            \small
            \node[rectangle, draw = black, dashed, minimum width = 0.75cm, minimum height = 0.75cm, line width=0.5mm] (r1) at (0.71,0.50) {};
            \node at (-0.25, 0.5){\includegraphics[keepaspectratio, width=0.5 in]{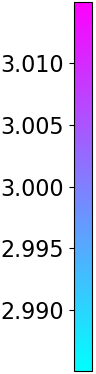}};
            \node at (-0.05, 0.5){\includegraphics[keepaspectratio, width=0.5 in]{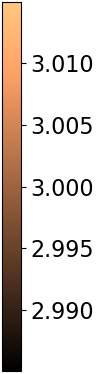}};
            \node at (0.51,0.5){\includegraphics[keepaspectratio, width=0.3 in]{"jupiter.png"}};
            \node at (1.21,0.5){\includegraphics[keepaspectratio, width=1.25 in]{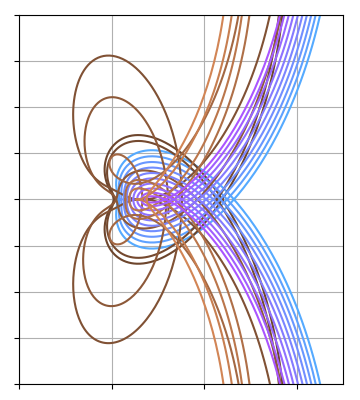}};
            \node[rectangle, draw = black, dashed, minimum width = 3.5cm, minimum height = 4.0cm, line width=0.5mm] (r2) at (1.21,0.5) {};
            \draw (r1.north west) -- (r2.north west);
            \draw[opacity=0.1] (r1.north east) -- (r2.north east);
            \draw (r1.south west) -- (r2.south west);
            \draw[opacity=0.1] (r1.south east) -- (r2.south east);
            \node[fill=white, opacity=1.0, text opacity=1, anchor=south] at (0.525,-0.0425) {$q_1$ [DU]};
            \node[fill=white, opacity=1.0, text opacity=1, anchor=south, rotate=90] at (0.125,0.5) {$q_2$ [DU]};
            \node[fill=white, opacity=1.0, text opacity=1, anchor=south] at (-0.145,-0.0425) {Jacobi Integral};
            \large
        \end{tikzonimage}
        \caption{A family of 3:4 and 5:6 resonant orbits possessing Jacobi integrals in the range [2.995, 3.005] is shown, with the orbits colored by the corresponding energy levels.}
        \label{fig: orbit by energy}
    \end{figure}
}
{
}

The Jacobi integral of a low-thrust trajectory changes whenever the spacecraft executes a maneuver, so we need to ensure that the invariant manifolds, to which we are comparing the trajectory puncture points, possess the same Jacobi integral. 
To do so, we first uniformly discretize the Jacobi integral interval between the initial orbit (3:4 resonant orbit, $C_{3:4}=2.995$) and the final orbit (5:6 resonant orbit, $C_{5:6}=3.005$) with steps of 0.001, and for each point in that interval, we compute the resonant orbit (and their invariant manifolds) possessing that Jacobi integral.
We store the information in a look-up table which will be used later in the subsequent analysis.
Ideally, we would like to compare every point on a trajectory with the pertinent resonant orbits at those energy levels, but doing so would be numerically intractable.
Fig. \ref{fig: orbit by energy} illustrates the resonant orbits used in this analysis.
To compute the resonant orbits, we use the initial conditions from the database developed by Restrepo and Russell \cite{restrepo_database_2018}.
For each resonant orbit shown, their invariant manifolds are also computed (not shown in the figure), and as before, approximately 10,000 puncture points are recorded for each invariant manifold.

\ifthenelse{\boolean{includefigures}}
{
    \begin{figure}[!htb]
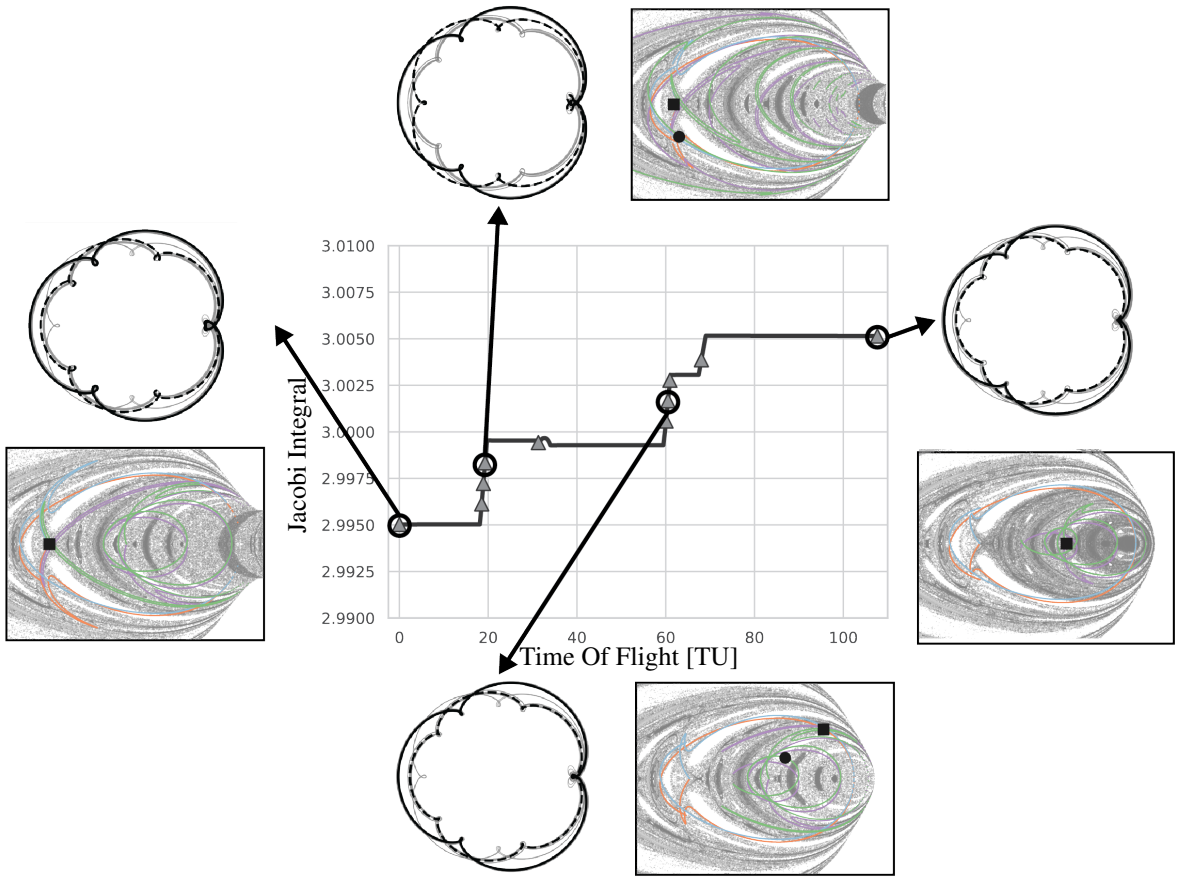

        \centering
        \begin{tikzonimage}[keepaspectratio, width=\linewidth]{"energy_orbit_section.png"}
        \node[fill=white, opacity=0.0, text opacity=1, anchor=south] at (0.533, 0.24) {Time Of Flight [TU]};
        \node[fill=white, opacity=0.0, text opacity=1, anchor=south, rotate=90] at (0.285,0.485) {Jacobi Integral};
        \end{tikzonimage}
        \caption{The Jacobi integral corresponding to an example non-robust solution is shown here, along with the trajectory points used for analysis highlighted with triangles.
        Each point possess the same Jacobi integral as a resonant orbit pair in our look-up table. 
        For example, the first point corresponds to a Jacobi integral of $\approx 2.995$.
        So, we compare its forward and backward integrated trajectory with the 3:4 and 5:6 resonant orbit invariant manifolds possessing the same energy.
        The 3:4 resonant orbit is shown with a solid line and the 5:6 resonant orbit with a dashed black line.
        The remaining orbits in the dictionary are shown in the background in \textcolor{Grey}{grey}.
        The corresponding Poincar\'{e} section $\mathcal{S}$ with the relevant puncture points is also shown below.}
        \label{fig: energy_orbit_section}
    \end{figure}
}
{
}

Figure \ref{fig: energy_orbit_section} illustrates the Jacobi integral associated with a representative non-robust solution. 
As the spacecraft trajectory evolves in time, its Jacobi integral undergoes a change every time the spacecraft executes a maneuver.
We begin by filtering points along the trajectory possessing a Jacobi integral close to that of one of the periodic orbits in our look-up table (in our case, we set this threshold to $1 \times 10^{-6}$).
From each filtered subset, we randomly select a point and integrate it forward and backward in time under natural dynamics until they intersect with the Poincar\'{e} section $\mathcal{S}$ to produce $\mathcal{P}(x)$ and $\mathcal{P}^{-1}(x)$ respectively.
We compare these puncture points to those of the invariant manifolds of the resonant orbits with the same Jacobi integral, ensuring that all points on a given Poincar\'{e} section possess the same energy.

\subsection{Distance Metrics on Poincar\'{e} Surface of Sections}
\label{subsection: analysis methods: distance metrics on poincare surface of sections}

\ifthenelse{\boolean{includefigures}}
{
    \begin{figure}[!htb]
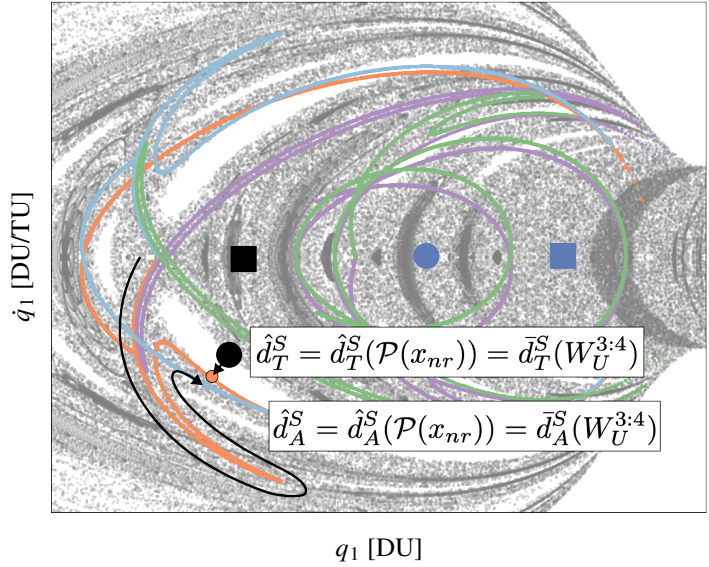

        \centering
        \begin{tikzonimage}[keepaspectratio, width=3.5 in]{"distance_metrics_example.png"}
        \node[fill=white, opacity=1.0, text opacity=1, anchor=south] at (0.5,-0.1) {$q_1$ [DU]};
        \node[fill=white, opacity=1.0, text opacity=1, anchor=south, rotate=90] at (0.0,0.5) {$\dot{q}_1$ [DU/TU]};
        \end{tikzonimage}
        \caption{An example Poincar\'{e} section is shown to visualize the distance metrics used in this study.
        In this particular frame, $\mathcal{P}(x_{nr})$ is closest to $W_{U}^{3:4}$.
        Because $\mathcal{P}^{-1}(x_{nr})$ is farther away from $W_{U}^{3:4}$, and $\mathcal{P}(x_{nr})$ is closest to $W_{U}^{3:4}$, $\bar{d}_T^\mathcal{S}(W_{U}^{3:4})$ is equal to $\hat{d}_T^\mathcal{S}(\mathcal{P}(x_{nr}))$.
        Furthermore, since all the other trajectory puncture points are farther away from \textit{any} of the invariant manifolds, this distance is also equal to $\hat{d}_T^\mathcal{S}$.
        Similar idea follows for the $\hat{d}_A^\mathcal{S}$.}
        \label{fig: distance metrics example}
    \end{figure}
}
{
}

By the process outlined above, it is possible to visualize the evolution of the trajectory puncture points on the Poincar\'{e} section $\mathcal{S}$, and compare them to the invariant manifolds of the pertinent resonant orbits.
However, in addition to qualitatively visualizing the solutions, we also aim to understand quantitatively if there is a difference in the behavior of the solutions with respect to the invariant manifolds, and in that regard, we introduce two distance metrics on $\mathcal{S}$ (see Fig. \ref{fig: distance metrics example} for a schematic).

\subsubsection{Orthogonal Distance To the Nearest Invariant Manifold $d_T^\mathcal{S}$}
\label{subsubsection: analysis methods: distance to the nearest invariant manifold puncture point}

First, we introduce a distance metric to quantify the degree to which a given trajectory leverages the underlying invariant manifolds in our problem.
We define $d_T^\mathcal{S}(x^{\mathcal{P}}, \mathcal{W})$ as the shortest orthogonal distance between $x^{\mathcal{P}}$, a trajectory puncture point, and $\mathcal{W}$, a set of invariant manifold puncture points on $\mathcal{S}$.
For example, $d_T^\mathcal{S}(\mathcal{P}(x_{nr}), W_{U}^{3:4})$ represents the shortest distance between the non-robust forward-integrated trajectory puncture point $\mathcal{P}(x_{nr})$ and the unstable manifold puncture points of the 3:4 resonant orbit $W_{U}^{3:4}$. 
We denote the shortest distance from a given trajectory puncture point $x^{\mathcal{P}}$ to \textit{any} of the invariant manifolds as: 
\begin{equation}
    \hat{d}_T^\mathcal{S}(x^{\mathcal{P}}) \equiv \min_{W_i^\gamma} d_T^\mathcal{S}(x^{\mathcal{P}}, W_i^\gamma)
    \label{eq:minimum distance to from a specific trajectory puncture point}
\end{equation}
where $x^{\mathcal{P}}$ may refer to either the forward-integrated trajectory point $\mathcal{P}(x)$ or the backward-integrated trajectory point $\mathcal{P}^{-1}(x)$.
Conversely, we denote the shortest distance from a given invariant manifold $W_i^\gamma$ of a resonant orbit to any of the trajectory puncture points as: 
\begin{equation}
    \bar{d}_T^\mathcal{S}(W_i^\gamma) \equiv \min_{x^{\mathcal{P}} \in \{\mathcal{P}(x), \mathcal{P}^{-1}(x)\}} d_T^\mathcal{S}(x^{\mathcal{P}}, W_i^\gamma)
    \label{eq:minimum distance to from any trajectory puncture point to a specific invariant manifold}
\end{equation}
where $W_i^\gamma$ may refer to $W_U^{3:4}$, $W_S^{3:4}$, $W_U^{5:6}$ or $W_S^{5:6}$.

\subsubsection{Distance Along the Nearest Invariant Manifold $d_A^\mathcal{S}$}
\label{subsubsection: analysis methods: distance along the nearest invariant manifold puncture point}

We introduce an additional distance metric to quantify the distance a trajectory must traverse along an invariant manifold to reach the nearest resonant orbit. 
Let $w \in \mathcal{W}$ represent a point in the set containing the invariant manifold puncture points $\mathcal{W}$. 
The arc length between $w$ and the separatrix, measured along the invariant manifold, quantifies the distance the spacecraft will have to coast to reach the corresponding resonant orbit. 
We denote this arc length as $d_A^\mathcal{S}(x^{\mathcal{P}}, w)$.
Of particular interest is the distance along the invariant manifold of an invariant manifold puncture point that is orthogonally closest to the trajectory puncture point, i.e., the invariant manifold puncture points that yields $\hat{d}_T^\mathcal{S}(x^{\mathcal{P}})$, which we denote as:
\begin{equation}
    \hat{d}_A^\mathcal{S}(x^{\mathcal{P}}) \equiv d_A^\mathcal{S}(x^{\mathcal{P}}, w) \ \text{s.t.} \ w = \underset{W_j^i}{\arg\min} \ d_T^\mathcal{S}(x^{\mathcal{P}}, W_j^i)
    \label{eq:minimum distance}
\end{equation}
Conversely, we denote the arc length along a given invariant manifold puncture point set $W_i^\gamma$ from the point in the invariant manifold puncture point set closest to the trajectory puncture points as: 
\begin{equation}
    \bar{d}_A^\mathcal{S}(W_i^\gamma) \equiv \min_{x^{\mathcal{P}} \in \{\mathcal{P}(x), \mathcal{P}^{-1}(x)\}} d_A^\mathcal{S}(x^{\mathcal{P}}, w) \ \text{s.t.} \ w = \underset{W_j^i}{\arg\min} \ d_T^\mathcal{S}(x^{\mathcal{P}}, W_j^i)
    \label{eq:minimum distance along from any trajectory puncture point to a specific invariant manifold}
\end{equation}

Finally, we introduce the following definitions to represent the minimum distance between the invariant manifold puncture points and the trajectory puncture points, irrespective of whether they were forward-integrated or backward-integrated in time:
\begin{equation}
    \hat{d}_T^\mathcal{S} \equiv \min \{ d_T^\mathcal{S}(\mathcal{P}(x)), d_T^\mathcal{S}(\mathcal{P}^{-1}(x)) \}
    \label{eq:minimum distance to from any trajectory puncture point}
\end{equation}
\begin{equation}
    \hat{d}_A^\mathcal{S} \equiv \min \{ d_A^\mathcal{S}(\mathcal{P}(x)), d_A^\mathcal{S}(\mathcal{P}^{-1}(x)) \}
    \label{eq:minimum distance along from any trajectory puncture point}
\end{equation}
The metric $\hat{d}_T^\mathcal{S}$ provides a measure of how far, at any given time, the trajectory is from \textit{any} of the invariant manifolds. 
By recording this distance at multiple time snapshots across different solutions, we can assess how effectively different solution families exploit the underlying invariant manifolds for minimum-fuel transfers, thereby providing insights into the role of dynamical structures in the optimization process.
On the other hand, the metric $\hat{d}_A^\mathcal{S}$ quantifies the distance the spacecraft will have to coast to reach the nearest resonant orbit at any given energy level.
This distance allows us to elucidate how fast the trajectory can coast along the nearest invariant manifold, and thereby, provide deeper insights into how these trajectories are leveraging these structures.

%% file: sections/results_and_discussion.tex
\section{Results and Discussion}
\label{section: results and discussion}

\ifthenelse{\boolean{includefigures}}
{
    \begin{figure}[!htb]
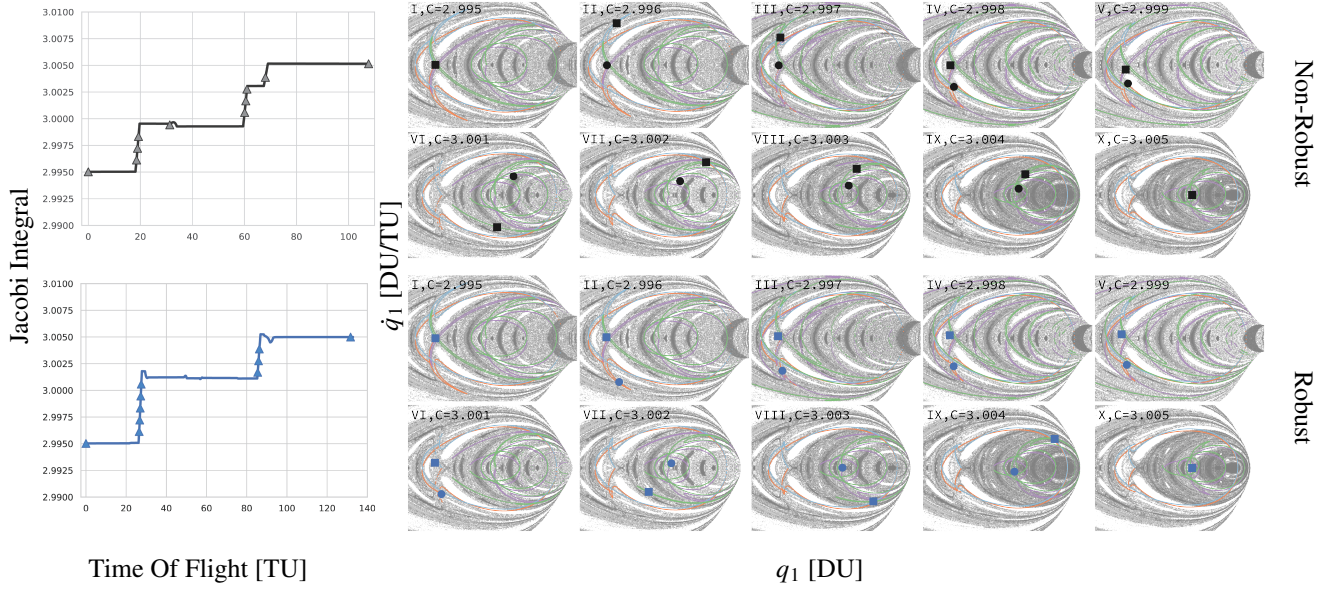

        \begin{tikzonimage}[keepaspectratio, width=\linewidth]{"frames_nonrobust_vs_robust.png"}
        \node[fill=white, opacity=0.0, text opacity=1, anchor=south] at (0.125,-0.075) {Time Of Flight [TU]};
        \node[fill=white, opacity=0.0, text opacity=1, anchor=south, rotate=90] at (0.0,0.5) {Jacobi Integral};
        \node[fill=white, opacity=0.0, text opacity=1, anchor=south] at (0.625,-0.075) {$q_1$ [DU]};
        \node[fill=white, opacity=0.0, text opacity=1, anchor=south, rotate=90] at (0.2975,0.5) {$\dot{q}_1$ [DU/TU]};
        \node[fill=white, opacity=0.0, text opacity=1, anchor=south, rotate=-90] at (1.0,0.75) {Non-Robust};
        \node[fill=white, opacity=0.0, text opacity=1, anchor=south, rotate=-90] at (1.0,0.25) {Robust};
        \end{tikzonimage}
        \caption{
        Representative snapshots in time demonstrating the temporal evolution of puncture points on $\mathcal{S}$ are shown. 
        Forward-integrated non-robust (robust) solution puncture points are denoted by black (\textcolor{TealBlue}{blue}) circles and backward-integrated non-robust (robust) solution puncture points by black (\textcolor{TealBlue}{blue}) squares.
        }
        \label{fig: frames_nonrobust_vs_robust}
    \end{figure}
}
{
}

\subsection{Qualitative Analysis}
\label{subsection: results and discussion: qualitative analysis}

To differentiate between robust and non-robust solutions, we analyze, both qualitatively and quantitatively, snapshots across the evolution of a specific solution, and compare the robust puncture points with their non-robust counterparts.
Intuitively, one might expect robust solutions to be less efficient than non-robust ones, leading to larger deviations from the underlying invariant manifolds, and this is certainly true for the solutions shown in Fig. \ref{fig: frames_nonrobust_vs_robust} which corresponds to an optimal non-robust solution and an optimal robust solution ($\tau_1$ coincides with the beginning of the 44-th control segment, and $\Delta \tau$ = 2.5 TU).
Looking at the temporal evolution of the Jacobi integral over time, we notice that the rate of change of the Jacobi integral of the non-robust solution is mostly non-negative, while that of the robust solution exhibits more frequent fluctuations.
A majority of robust solutions display similar behavior, and this observation aligns with our hypothesis that robust solutions can sometime undergo inefficient maneuvers to attain feasibility.

As seen in Fig. \ref{fig: frames_nonrobust_vs_robust}, the robust trajectory puncture points follow a different path than the non-robust puncture points but both seem to be relatively well-aligned with the invariant manifold puncture points.
We now analyze the non-robust solution in more detail.
As one would expect, both $\mathcal{P}(x_{nr})$ and $\mathcal{P}^{-1}(x_{nr})$ begin at the 3:4 separatrix (Frame I).
$\mathcal{P}(x_{nr})$ remains close to this point (Frames II-III), before transitioning to $W_U^{3:4}$ which it flows along (Frames IV-V).
Then, it transfers between invariant manifolds (Frames VI-IX), maintaining close proximity to both $W_U^{5:6}$ and $W_S^{5:6}$, until it eventually reaches the 5:6 separatrix (Frame X).
On the other hand, $\mathcal{P}^{-1}(x_{nr})$ immediately latches onto $W_S^{3:4}$ which it flows along (Frames II-III), before it returns back to the separatrix (Frame IV-V).
Then, it transitions to $W_S^{5:6}$ which it flows along (Frames VI-IX) until it reaches the 5:6 separatrix (Frame X).
Between Frames I-V, the spacecraft has optimized its thrust to move $\mathcal{P}(x_{nr})$ along $W_U^{3:4}$ while aligning $\mathcal{P}^{-1}(x_{nr})$ for the subsequent transition into $W_S^{5:6}$.
Then, between Frames VI-X, the thrust is used to transition $\mathcal{P}^{-1}(x_{nr})$ along $W_S^{5:6}$, while transferring $\mathcal{P}(x_{nr})$ between different invariant manifolds until they settle at the 5:6 separatrix.

Moving on to the robust solution, naturally, both $\mathcal{P}(x_{r})$ and $\mathcal{P}^{-1}(x_{r})$ naturally begin at the 3:4 separatrix (Frame I).
$\mathcal{P}(x_{r})$ aligns itself with the $W_U^{3:4}$ which it flows along (Frames II-VI). 
In the subsequent frames, it `oscillates' about the 5:6 separatrix (Frames VII-IX) before eventually setting at the separatrix (Frame X).
$\mathcal{P}^{-1}(x_{r})$, on the other hand, follows a rather interesting path.
We do not see any major change during the initial frames (Frames I-VI).
Then, we see it abruptly transition to $W_S^{5:6}$ (Frame VII) which it flows along (Frames VII-IX) until it reaches the 5:6 separatrix (Frame X).
Between Frames I-VI, the spacecraft has optimized its thrust to first move $\mathcal{P}(x_{r})$ along $W_U^{3:4}$ while making minimal changes to $\mathcal{P}^{-1}(x_{r})$, and then, in Frame VII, to move both $\mathcal{P}(x_{r})$ and $\mathcal{P}^{-1}(x_{r})$ over to $W_S^{5:6}$.
Finally, between Frames VII-X, the thrust is mostly used to move $\mathcal{P}(x_{r})$ and $\mathcal{P}^{-1}(x_{r})$ along $W_S^{5:6}$ (Frames VII-IX) until they settle at the 5:6 separatrix.

\subsection{Quantitative Analysis}
\label{subsection: results and discussion: quantitative analysis}

From the analysis of this solution, we observe that both robust and non-robust trajectories exhibit certain qualitative similarities and differences. 
Notably, the robust solution appears to flow along the invariant manifolds whereas the non-robust solution appears to transfer between the invariant manifolds (e.g., Frames VI-IX). 
However, generalizing the relationship between trajectory and invariant manifold puncture points based solely on visual inspection of a single solution pair is difficult.
The observed trend in Fig. \ref{fig: frames_nonrobust_vs_robust} may be specific to this particular pair and may not be representative of the entire solution family.
Further, it is also important to note that the snapshots are non-uniform in time.
As a result, we often notice large `jumps' in the path of the trajectory puncture points (e.g., Frames V \& VI for the non-robust solution; Frames VI \& VII for the robust solution).
So, instead, by analyzing a collection of solutions, we hope to be able to `average out' the non-uniformity in time to globally characterize the solution trends in relation to the invariant manifolds.

\begin{table}[hbt!]
    \caption{Number of Solutions and Punctures}
    \label{tab: number_of_solutions_and_punctures}
    \centering
    \begin{tabular}{llllcccllll}
    \cline{1-7}
                                 & \multicolumn{2}{c}{\textbf{Non-Robust}}              & \multicolumn{1}{c}{\textbf{}}       & \multicolumn{1}{l}{\textbf{}}             & \multicolumn{2}{c}{\textbf{Robust}}                        &  &  &  &  \\ \cline{1-7}
                                 & Feasible                  & Optimal                  &                                     & \multicolumn{1}{l}{}                      & \multicolumn{1}{l}{Feasible} & \multicolumn{1}{l}{Optimal} &  &  &  &  \\ \cline{1-7}
    \textbf{Number of Solutions} & \multicolumn{1}{c}{21455} & \multicolumn{1}{c}{9393} &                                     &                                           & 6789                         & 958                         &  &  &  &  \\ \cline{1-7}
                                 & \multicolumn{1}{c}{}      & \multicolumn{1}{c}{}     & \multirow{7}{*}{$\Delta \tau$ (TU)} & 0.5                                       & 1051                         & 171                         &  &  &  &  \\
                                 & \multicolumn{1}{c}{}      & \multicolumn{1}{c}{}     &                                     & 1.0                                       & 932                          & 140                         &  &  &  &  \\
                                 & \multicolumn{1}{c}{}      & \multicolumn{1}{c}{}     &                                     & 2.5                                       & 1016                         & 156                         &  &  &  &  \\
                                 & \multicolumn{1}{c}{}      & \multicolumn{1}{c}{}     &                                     & 5.0                                       & 1002                         & 153                         &  &  &  &  \\
                                 & \multicolumn{1}{c}{}      & \multicolumn{1}{c}{}     &                                     & 10.0                                      & 969                          & 133                         &  &  &  &  \\
                                 & \multicolumn{1}{c}{}      & \multicolumn{1}{c}{}     &                                     & 15.0                                      & 951                          & 121                         &  &  &  &  \\
                                 & \multicolumn{1}{c}{}      & \multicolumn{1}{c}{}     &                                     & 30.0                                      & 868                          & 84                          &  &  &  &  \\ \cline{4-7}
                                 &                           &                          & \multirow{2}{*}{$\tau_1$}             & \multicolumn{1}{l}{Forward Shooting Arc}  & 2850                         & 224                         &  &  &  &  \\
                                 &                           &                          &                                     & \multicolumn{1}{l}{Backward Shooting Arc} & 3939                         & 734                         &  &  &  &  \\ \cline{1-7}
    \textbf{Number of Punctures} & 168722                    & 143890                   &                                     & \multicolumn{1}{l}{}                      & 65895                        & 17202                       &  &  &  &  \\ \cline{1-7}
                                 &                           &                          & \multirow{7}{*}{$\Delta \tau$ (TU)} & 0.5                                       & 10068                        & 3071                        &  &  &  &  \\
                                 &                           &                          &                                     & 1.0                                       & 8902                         & 2399                        &  &  &  &  \\
                                 &                           &                          &                                     & 2.5                                       & 9918                         & 2914                        &  &  &  &  \\
                                 &                           &                          &                                     & 5.0                                       & 9702                         & 2753                        &  &  &  &  \\
                                 &                           &                          &                                     & 10.0                                      & 9459                         & 2433                        &  &  &  &  \\
                                 &                           &                          &                                     & 15.0                                      & 9275                         & 2104                        &  &  &  &  \\
                                 &                           &                          &                                     & 30.0                                      & 8571                         & 1501                        &  &  &  &  \\ \cline{4-7}
                                 &                           &                          & \multirow{2}{*}{$\tau_1$}             & \multicolumn{1}{l}{Forward Shooting Arc}  & 33014                        & 4300                        &  &  &  &  \\
                                 &                           &                          &                                     & \multicolumn{1}{l}{Backward Shooting Arc} & 32881                        & 12902                       &  &  &  &  \\ \cline{1-7}
    \end{tabular}
\end{table}

To assess whether the observed trend persists across solution families, we first generate a set of feasible and optimal solutions for both the non-robust and robust problems, varying the parameters $\tau_1$ and $\Delta \tau$. 
The optimal solutions, a subset of the feasible solutions, comprises the solutions that meet the optimizer's optimality criteria, specifically those exiting with \texttt{SNOPT Exit Info 1} \cite{gill_snopt_usermanual_2018}. 
The number of solutions and puncture points used in this study are presented in Table \ref{tab: number_of_solutions_and_punctures}.

In this section, we provide a detailed quantitative analysis to complement the qualitative observations.
Using the distance metrics discussed in \S \ref{subsection: analysis methods: distance metrics on poincare surface of sections}, we aim to examine the differences between the non-robust and robust solutions relative to the invariant manifolds.
We begin by computing $\hat{d}_T^\mathcal{S}$ and $\hat{d}_A^\mathcal{S}$ for the non-robust and the robust solutions, and then compare the statistics between the feasible solutions and the optimal solutions in each category (Figs. \ref{fig: nonrob_distance}, \ref{fig: rob_distance}).

\ifthenelse{\boolean{includefigures}}
{
    \begin{figure}[!htb]
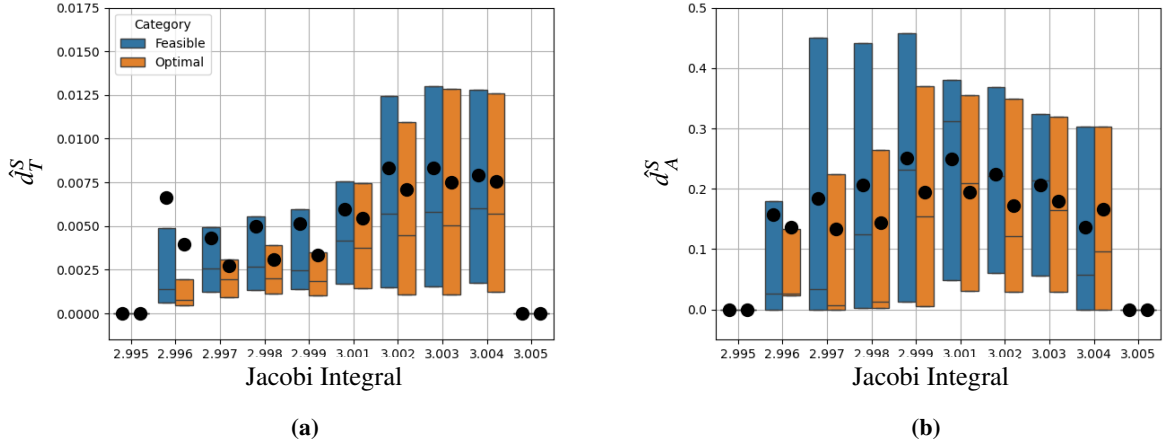

        \centering
        \begin{subfigure}[b]{0.5\textwidth}
            \centering
            \begin{tikzonimage}[keepaspectratio, width=3 in]{"nonrob_distance_to.png"}
            \node[fill=white, opacity=1.0, text opacity=1, anchor=south] at (0.5,-0.03) {Jacobi Integral};
            \node[fill=white, opacity=1.0, text opacity=1, anchor=south, rotate=90] at (0.025,0.5) {$\hat{d}^S_T$};
            \end{tikzonimage}
            \caption{}
            \label{fig: nonrob_distance_to}
        \end{subfigure}%
        \begin{subfigure}[b]{0.5\textwidth}
            \centering
            \begin{tikzonimage}[keepaspectratio, width=3 in]{"nonrob_distance_along.png"}
            \node[fill=white, opacity=1.0, text opacity=1, anchor=south] at (0.5,-0.03) {Jacobi Integral};
            \node[fill=white, opacity=1.0, text opacity=1, anchor=south, rotate=90] at (0.075,0.5) {$\hat{d}^S_A$};
            \end{tikzonimage}
            \caption{}
            \label{fig: nonrob_distance_along}
        \end{subfigure}
        \caption{$\hat{d}_T^\mathcal{S}$ and $\hat{d}_A^\mathcal{S}$ for the non-robust solution family are shown with the overall feasible solution set denoted by \textcolor{python_deep_blue}{blue}, and the optimal solution subset denoted by \textcolor{python_deep_orange}{orange}.
        The circles represent the mean, the horizontal lines represent the median, and the bars represent the interquartile range.}
        \label{fig: nonrob_distance}
    \end{figure}
}
{
}

\ifthenelse{\boolean{includefigures}}
{
    \begin{figure}[!htb]
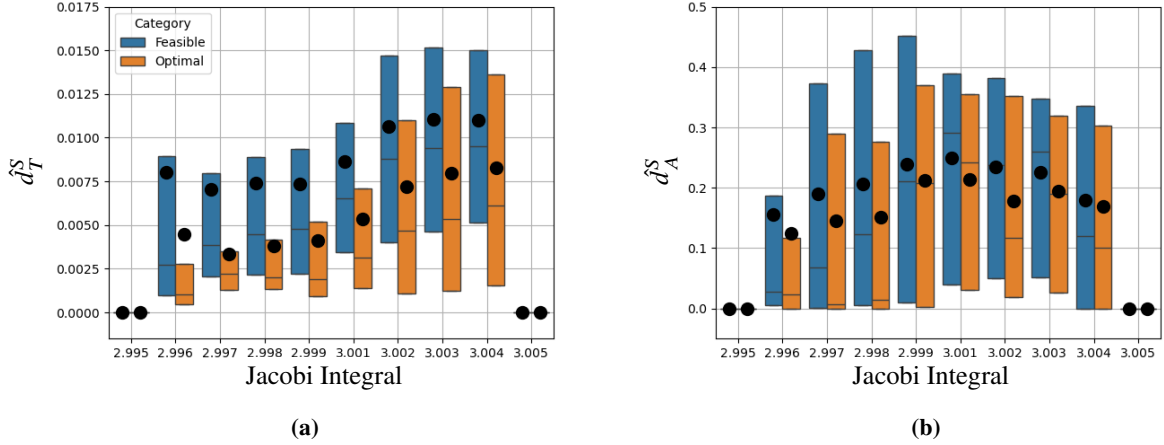

        \centering
        \begin{subfigure}[b]{0.5\textwidth}
            \centering
            \begin{tikzonimage}[keepaspectratio, width=3 in]{"rob_distance_to.png"}
            \node[fill=white, opacity=1.0, text opacity=1, anchor=south] at (0.5,-0.03) {Jacobi Integral};
            \node[fill=white, opacity=1.0, text opacity=1, anchor=south, rotate=90] at (0.025,0.5) {$\hat{d}^S_T$};
            \end{tikzonimage}
            \caption{}
            \label{fig: rob_distance_to}
        \end{subfigure}%
        \begin{subfigure}[b]{0.5\textwidth}
            \centering
            \begin{tikzonimage}[keepaspectratio, width=3 in]{"rob_distance_along.png"}
            \node[fill=white, opacity=1.0, text opacity=1, anchor=south] at (0.5,-0.03) {Jacobi Integral};
            \node[fill=white, opacity=1.0, text opacity=1, anchor=south, rotate=90] at (0.075,0.5) {$\hat{d}^S_A$};
            \end{tikzonimage}
            \caption{}
            \label{fig: rob_distance_along}
        \end{subfigure}
        \caption{
        $\hat{d}_T^\mathcal{S}$ and $\hat{d}_A^\mathcal{S}$ for the robust solution family are shown, which contain solutions with varying $\Delta \tau$ and $\tau_1$.
        }
        \label{fig: rob_distance}
    \end{figure}
}
{
}

As expected, $\hat{d}_T^\mathcal{S}$ is zero at the first and the last frames, since the trajectory puncture points at those frames coincide with the separatrices of the 3:4 and 5:6 resonant orbits respectively, and therefore the corresponding $\hat{d}_T^\mathcal{S}$ should be exactly zero.
As we progress through the energy levels, the mean $\hat{d}_T^\mathcal{S}$ increases for both solution categories.
Since the transfer is from a 3:4 resonant orbit farther from Europa to a 5:6 resonant orbit closer to Europa, we anticipate greater sensitivity in the solutions during the later energy levels, which correspond to points along the trajectory closer to Europa, and therefore subject to dynamics highly sensitive to perturbations (Fig. \ref{fig: orbit by energy}).
Accordingly, we expect the puncture points corresponding to the later energy levels to show less reliance on the invariant manifolds compared to those at earlier levels, as the dynamics may be too chaotic for the optimizer to leverage the invariant manifolds effectively for a finite horizon minimum-fuel transfer.
This expectation is confirmed by the higher mean $\hat{d}_T^\mathcal{S}$ observed at the later energy levels relative to the earlier ones.
Throughout the analysis, however, the optimal solutions consistently remain closer to the invariant manifolds on average compared to the feasible solutions, a trend that persists across energy levels and solution categories, suggesting that closer alignment with the invariant manifolds is necessary to achieve optimality.
Robust feasible solutions exhibit a slightly higher $\hat{d}_T^\mathcal{S}$ compared to the non-robust solutions, while the optimal set of robust solutions demonstrates a substantial decrease in $\hat{d}_T^\mathcal{S}$, bringing it nearly in line with those for the non-robust solutions.

The statistics for $\hat{d}_A^\mathcal{S}$ show a similar trend between the feasible and optimal solutions, for both the non-robust and the robust solution categories.
For both solution categories, the $\hat{d}_A^\mathcal{S}$ remain relatively small during the initial energy levels, which increases as we progress through the energy levels, before decreasing again at the later energy levels.
Given that the solutions involve transfers between resonant orbits, we expect $\hat{d}_A^\mathcal{S}$ to remain small when the trajectory is in proximity to the resonant orbits, which justifies the lower $\hat{d}_A^\mathcal{S}$ in the earlier and the later energy levels when the spacecraft is departing the 3:4 resonant orbit and entering the 5:6 resonant orbit respectively, and higher in between.

To compare the robust solutions with the non-robust solutions as we vary $\tau_1$ and $\Delta \tau$, we calculate the fold change in the mean distance metrics for the robust solutions relative to their non-robust counterparts. 
To visualize this change, we use the $\text{log}_2$ of the fold change as a metric, where a value of zero indicates that the distance metrics for the robust solutions and the non-robust solutions are the same, a positive value indicates an increase and a negative value indicates a decrease.
We remove the first and energy levels from the subsequent analysis since they correspond to the separatrices of the 3:4 and 5:6 resonant orbits respectively.
We begin by analyzing the solutions grouped according to $\Delta \tau$, noting that each group contains solutions with varying $\tau_1$. 
Subsequently, we analyze the solutions by grouping them according to $\tau_1$, where each group similarly includes solutions with different $\Delta \tau$.

\subsubsection{Dependence on $\Delta \tau$}
\label{subsubsection: results and discussion: quantitative analysis: dependence on delta tau}

\ifthenelse{\boolean{includefigures}}
{
    \begin{figure}[!htb]
        \centering
        \begin{subfigure}[b]{0.5\textwidth}
            \centering
            \begin{tikzonimage}[keepaspectratio, width=3 in]{"mte_duration_distance_to_feasible.png"}
            \node[fill=white,opacity=1.0] at (0.335,0.78){\includegraphics[keepaspectratio, width=1.125 in]{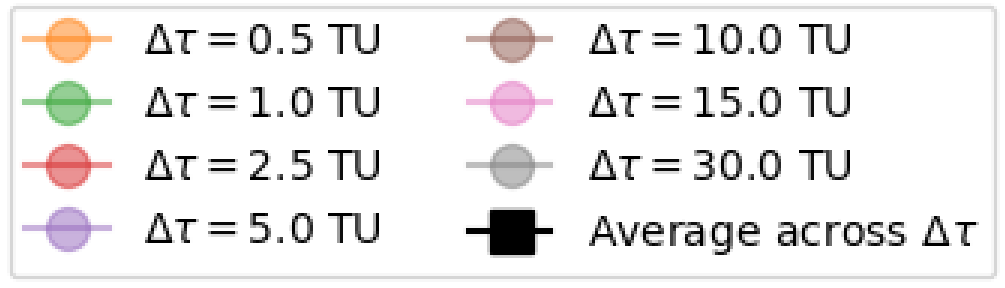}};
            \node[fill=white, opacity=1.0, text opacity=1, anchor=south] at (0.5,-0.03) {Jacobi Integral};
            \node[fill=white, opacity=1.0, text opacity=1, anchor=south, rotate=90] at (0.075,0.5) {$\text{log}_2 \left(\frac{\expectation{\omega}{\hat{d}^S_T (x_r)}}{\expectation{\omega}{\hat{d}^S_T (x_{nr})}}\right)$};
            \end{tikzonimage}
            \caption{Feasible Solutions}
            \label{fig: mte_duration_distance_to_feasible}
        \end{subfigure}%
        \begin{subfigure}[b]{0.5\textwidth}
            \centering
            \begin{tikzonimage}[keepaspectratio, width=3 in]{"mte_duration_distance_to_optimal.png"}
            \node[fill=white, opacity=1.0, text opacity=1, anchor=south] at (0.5,-0.03) {Jacobi Integral};
            \node[fill=white, opacity=1.0, text opacity=1, anchor=south, rotate=90] at (0.075,0.5) {$\text{log}_2 \left(\frac{\expectation{\omega}{\hat{d}^S_T (x_r)}}{\expectation{\omega}{\hat{d}^S_T (x_{nr})}}\right)$};
            \end{tikzonimage}
            \caption{Optimal Solutions}
            \label{fig: mte_duration_distance_to_optimal}
        \end{subfigure}
        \caption{
        Fold change in $\hat{d}_T^\mathcal{S}$ for robust solutions with varying $\Delta \tau$ relative to non-robust solutions
        }
        \label{fig: mte_duration_distance_to}
    \end{figure}
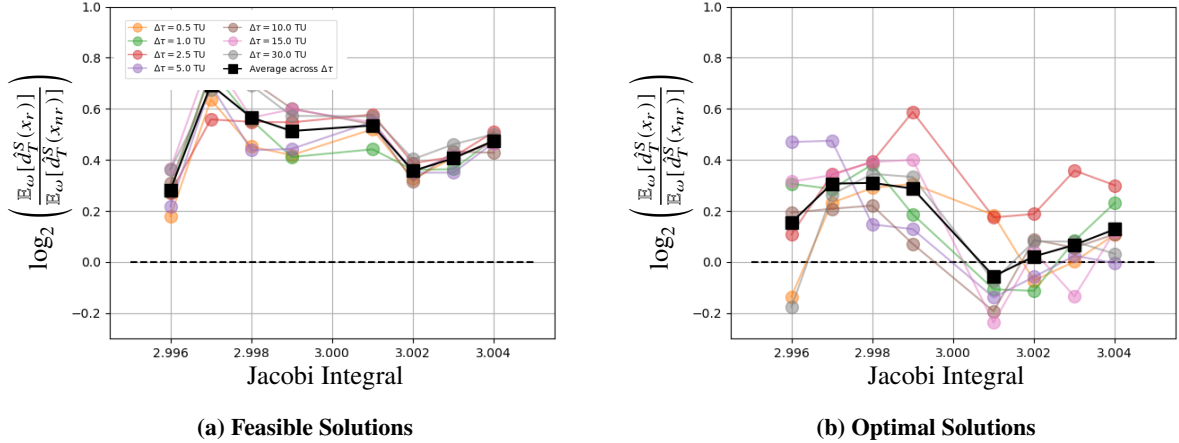
}
{
}

In this section, we explore how the relationship to the invariant manifolds changes as we vary $\Delta \tau$ for the robust solutions. 
If we consider the entire feasible solution set, the mean $\hat{d}_T^\mathcal{S}$ across $\Delta \tau$ for the robust solutions is greater than that for the non-robust solutions (Fig. \ref{fig: mte_duration_distance_to_feasible}). 
However, if we focus on the optimal solution subset, the mean $\hat{d}_T^\mathcal{S}$ across $\Delta \tau$ for the robust solutions significantly diminishes, and appears to resemble the non-robust solutions more closely, especially toward the later energy levels (Fig. \ref{fig: mte_duration_distance_to_optimal}). 
We even notice some robust solutions categories exhibiting closer alignment to the invariant manifolds than the non-robust solutions.
This observation is particularly significant because, as we have seen before, the solutions rely less on invariant manifolds at the later energy levels. 
The fact that robust solutions exhibit comparable, and in some cases stronger, alignment with these invariant manifolds at the later energy levels, relative to the non-robust solutions, suggests an important insight. 
Despite the chaotic dynamics, robust optimal solutions are still effectively leveraging the invariant manifolds, almost as closely as the non-robust solutions.
The results suggest that while feasible robust solutions may generally deviate from the invariant manifolds, the optimal ones tend to shadow the invariant manifolds almost as closely as the non-robust optimal solutions.

\ifthenelse{\boolean{includefigures}}
{
    \begin{figure}[!htb]
        \centering
        \begin{subfigure}[b]{0.5\textwidth}
            \centering
            \begin{tikzonimage}[keepaspectratio, width=3 in]{"mte_duration_distance_along_feasible.png"}
            \node[fill=white,opacity=1.0] at (0.335,0.78){\includegraphics[keepaspectratio, width=1.125 in]{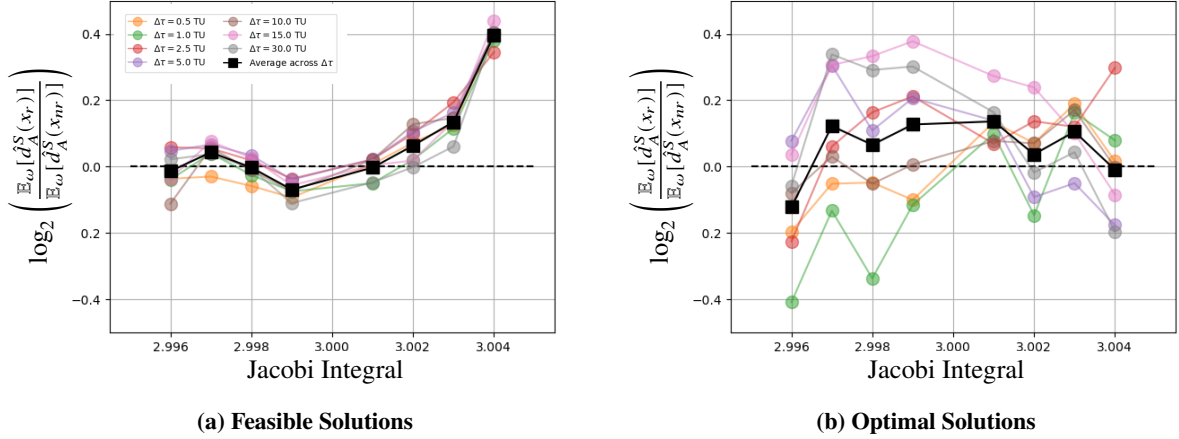}};
            \node[fill=white, opacity=1.0, text opacity=1, anchor=south] at (0.5,-0.03) {Jacobi Integral};
            \node[fill=white, opacity=1.0, text opacity=1, anchor=south, rotate=90] at (0.075,0.5) {$\text{log}_2 \left(\frac{\expectation{\omega}{\hat{d}^S_A (x_r)}}{\expectation{\omega}{\hat{d}^S_A (x_{nr})}}\right)$};
            \end{tikzonimage}
            \caption{Feasible Solutions}
            \label{fig: mte_duration_distance_along_feasible}
        \end{subfigure}%
        \begin{subfigure}[b]{0.5\textwidth}
            \centering
            \begin{tikzonimage}[keepaspectratio, width=3 in]{"mte_duration_distance_along_optimal.png"}
            \node[fill=white, opacity=1.0, text opacity=1, anchor=south] at (0.5,-0.03) {Jacobi Integral};
            \node[fill=white, opacity=1.0, text opacity=1, anchor=south, rotate=90] at (0.075,0.5) {$\text{log}_2 \left(\frac{\expectation{\omega}{\hat{d}^S_A (x_r)}}{\expectation{\omega}{\hat{d}^S_A (x_{nr})}}\right)$};
            \end{tikzonimage}
            \caption{Optimal Solutions}
            \label{fig: mte_duration_distance_along_optimal}
        \end{subfigure}
        \caption{
        Fold change in $\hat{d}_A^\mathcal{S}$ for robust solutions with varying $\Delta \tau$ relative to non-robust solutions
        }
        \label{fig: mte_duration_distance_along}
    \end{figure}
}
{
}

An interesting trend emerges in the second distance metric. During the initial energy levels, the robust optimal solutions exhibit a larger mean $\hat{d}_A^\mathcal{S}$ across $\Delta \tau$ relative to the non-robust optimal solutions (Fig. \ref{fig: mte_duration_distance_along_optimal}), particularly when compared to the overall feasible solution set (Fig. \ref{fig: mte_duration_distance_along_feasible}). 
At the initial energy levels, the mean $\hat{d}_A^\mathcal{S}$ across $\Delta \tau$ for the robust and non-robust feasible solutions are similar. 
However, as we progress through the energy levels, the ratio decreases slightly before increasing.
If we focus only on the optimal solutions subset, we find that this distance metric for the robust solutions is lower at the first and last frames but shows a modest increase in between. 
As explained before, this is the behavior we expect since solutions are typically closer to resonant orbits during departure from the 3:4 resonant orbit (i.e., the early frames) and arrival at the 5:6 resonant orbit (i.e., the later frames).
The behavior of the robust optimal solutions observed in this distance metric suggests stronger alignment with this hypothesis.

\ifthenelse{\boolean{includefigures}}
{
    \begin{figure}[!htb]
        \centering
        \begin{subfigure}[b]{0.5\textwidth}
            \centering
            \begin{tikzonimage}[keepaspectratio, width=3 in]{"rob_vs_nonrob_distance_to_optimal_unstable_duration.png"}
            \node[fill=white,opacity=0.85] at (0.335,0.775){\includegraphics[keepaspectratio, width=1.125 in]{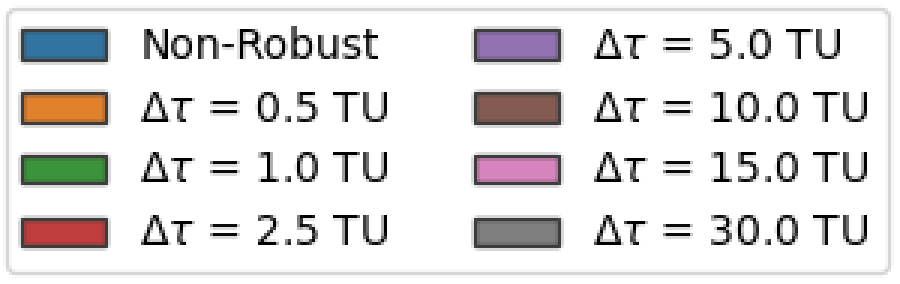}};
            \node[fill=white, opacity=1.0, text opacity=1, anchor=south] at (0.5,-0.03) {Jacobi Integral};
            \node[fill=white, opacity=1.0, text opacity=1, anchor=south, rotate=90] at (0.05,0.5) {$\bar{d}_T^S(W_U^{3:4})$};
            \end{tikzonimage}
            \caption{}
            \label{fig: rob_vs_nonrob_distance_to_optimal_unstable_duration}
        \end{subfigure}%
        \begin{subfigure}[b]{0.5\textwidth}
            \centering
            \begin{tikzonimage}[keepaspectratio, width=3 in]{"rob_vs_nonrob_distance_to_optimal_stable_duration.png"}
            \node[fill=white, opacity=1.0, text opacity=1, anchor=south] at (0.5,-0.03)  {Jacobi Integral};
            \node[fill=white, opacity=1.0, text opacity=1, anchor=south, rotate=90] at (0.05,0.5) {$\bar{d}_T^S(W_S^{5:6})$};
            \end{tikzonimage}
            \caption{}
            \label{fig: rob_vs_nonrob_distance_to_optimal_stable_duration}
        \end{subfigure}
        \caption{
        $\bar{d}_T^S(W_U^{3:4})$ and $\bar{d}_T^S(W_S^{5:6})$ for optimal non-robust solutions and robust solutions with varying $\Delta \tau$
        }
        \label{fig: rob_vs_nonrob_distance_to_optimal_duration}
    \end{figure}
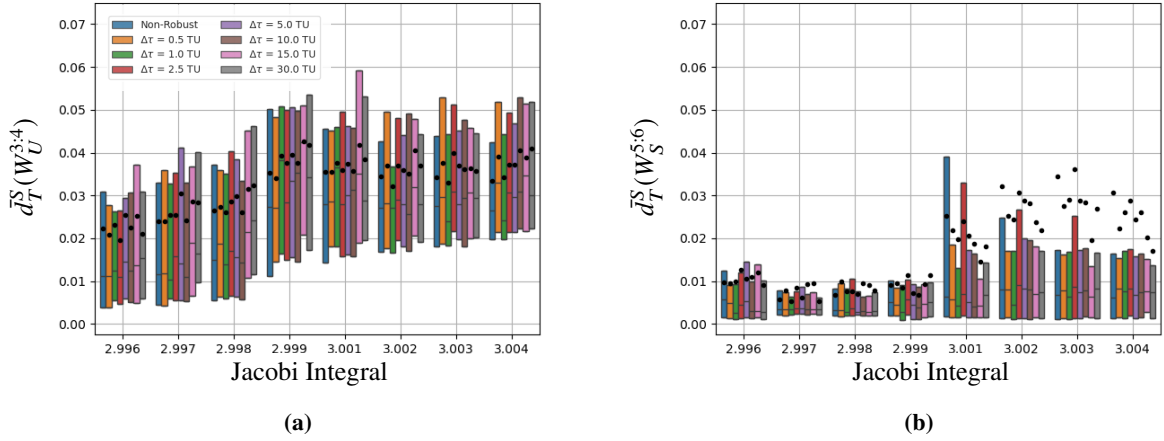
}
{
}

Looking at the distance metrics associated with particular invariant manifolds can also reveal useful insights into the solutions.
It can be intuitively surmised that the optimal low-thrust solutions for this problem shall flow along $W_U^{3:4}$ before eventually transitioning into $W_S^{5:6}$.
To test this hypothesis, we evaluate the orthogonal distance to $W_U^{3:4}$ and to $W_S^{5:6}$, i.e., $\bar{d}_T^S(W_U^{3:4})$ and $\bar{d}_T^S(W_S^{5:6})$ respectively, and the corresponding distance along the invariant manifold $\bar{d}_A^S(W_U^{3:4})$ and $\bar{d}_A^S(W_S^{5:6})$ for the optimal solutions in each category.
As indicated by the increasing mean $\bar{d}_T^S(W_U^{3:4})$ across energy levels, it is immediately obvious that the optimal solutions in all categories depart from $W_U^{3:4}$ as they traverse the energy levels (Fig. \ref{fig: rob_vs_nonrob_distance_to_optimal_unstable_duration}).
Although we do not observe a converse trend in $W_S^{5:6}$, it is important to note that the solutions remain consistently close to $W_S^{5:6}$ allowing them to leverage $W_S^{5:6}$ whenever necessary (Fig. \ref{fig: rob_vs_nonrob_distance_to_optimal_stable_duration}).

\ifthenelse{\boolean{includefigures}}
{
    \begin{figure}[!htb]
        \centering
        \begin{subfigure}[b]{0.5\textwidth}
            \centering
            \begin{tikzonimage}[keepaspectratio, width=3 in]{"rob_vs_nonrob_distance_along_optimal_unstable_duration.png"}
            \node[fill=white,opacity=0.85] at (0.335,0.775){\includegraphics[keepaspectratio, width=1.125 in]{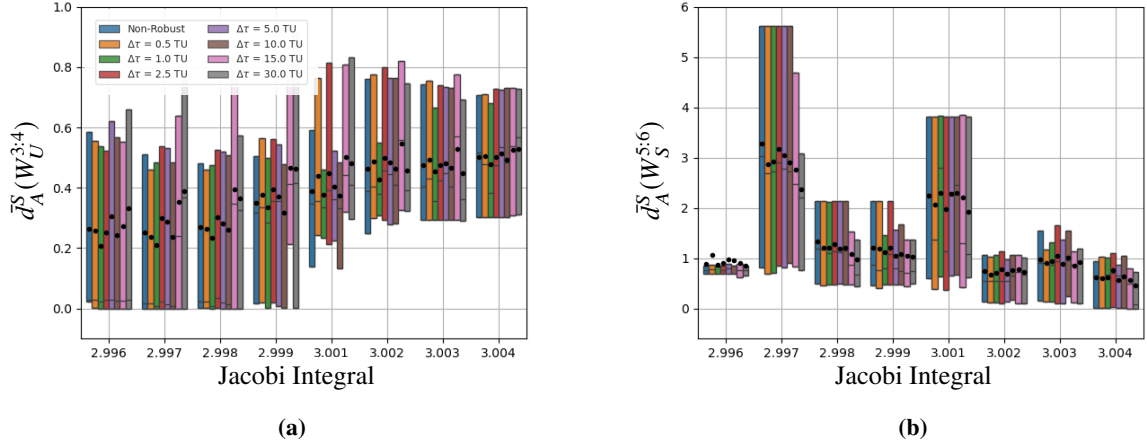}};
            \node[fill=white, opacity=1.0, text opacity=1, anchor=south] at (0.5,-0.03)  {Jacobi Integral};
            \node[fill=white, opacity=1.0, text opacity=1, anchor=south, rotate=90] at (0.075,0.5) {$\bar{d}_A^S(W_U^{3:4})$};
            \end{tikzonimage}
            \caption{}
            \label{fig: rob_vs_nonrob_distance_along_optimal_unstable_duration}
        \end{subfigure}%
        \begin{subfigure}[b]{0.5\textwidth}
            \centering
            \begin{tikzonimage}[keepaspectratio, width=3 in]{"rob_vs_nonrob_distance_along_optimal_stable_duration.png"}
            \node[fill=white, opacity=1.0, text opacity=1, anchor=south] at (0.5,-0.03)  {Jacobi Integral};
            \node[fill=white, opacity=1.0, text opacity=1, anchor=south, rotate=90] at (0.09,0.5) {$\bar{d}_A^S (W_S^{5:6})$};
            \end{tikzonimage}
            \caption{}
            \label{fig: rob_vs_nonrob_distance_along_optimal_stable_duration}
        \end{subfigure}
        \caption{
        $\bar{d}_A^S(W_U^{3:4})$ and $\bar{d}_A^S (W_S^{5:6})$ for optimal non-robust solutions and robust solutions with varying $\Delta \tau$
        }
        \label{fig: rob_vs_nonrob_distance_along_optimal_duration}
    \end{figure}
}
{
}

On the other hand, the trend in $\bar{d}_A^\mathcal{S}$ indicates that, on average, the solutions remain relatively close to the 3:4 resonant orbit as they traverse the initial energy levels. 
However, as the energy levels increase, the solutions begin to drift away from the 3:4 resonant orbit, reflecting a gradual departure from its influence (Fig. \ref{fig: rob_vs_nonrob_distance_along_optimal_unstable_duration}).
Conversely, as the solutions progress through higher energy levels, there is a noticeable shift in alignment toward the 5:6 resonant orbit. 
This suggests that the trajectories increasingly rely on the dynamical structures associated with the 5:6 resonant orbit as they approach the later stages of the transfer, which, once again, aligns with our expectations (Fig. \ref{fig: rob_vs_nonrob_distance_along_optimal_stable_duration}).

\ifthenelse{\boolean{includefigures}}
{
    \begin{figure}[!htb]
        \centering
        \begin{subfigure}[b]{0.5\textwidth}            
            \centering
            \begin{tikzonimage}[keepaspectratio, width=3 in]{"mte_duration_distance_to_along_feasible.png"}
            \node[fill=white,opacity=0.85] at (0.335,0.215){\includegraphics[keepaspectratio, width=1.125 in]{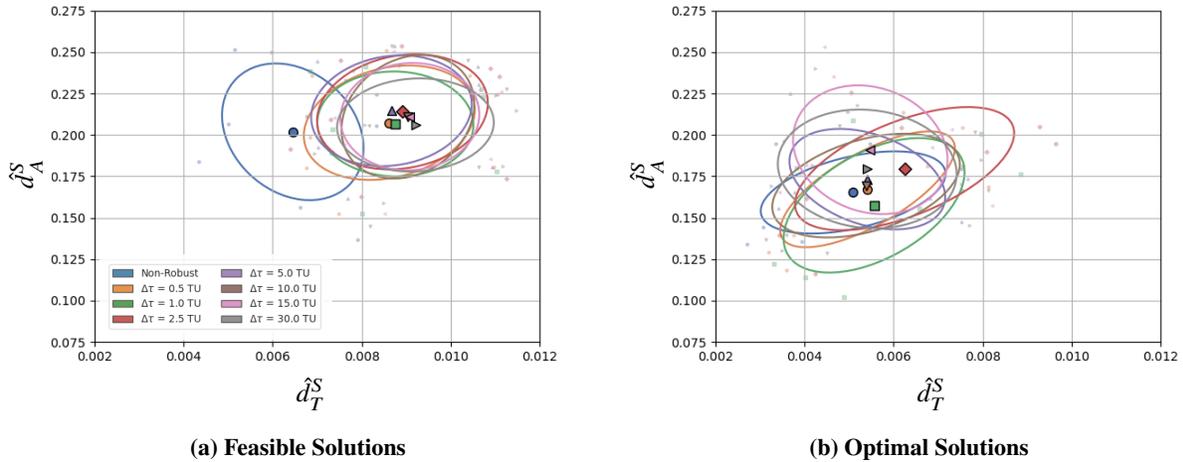}};
            \node[fill=white, opacity=1.0, text opacity=1, anchor=south] at (0.5,-0.07) {$\hat{d}^S_T$};
            \node[fill=white, opacity=1.0, text opacity=1, anchor=south, rotate=90] at (0.05,0.5) {$\hat{d}^S_A$};
            \end{tikzonimage}
            \caption{Feasible Solutions}
            \label{fig: mte_duration_distance_to_along_feasible}
        \end{subfigure}%
        \begin{subfigure}[b]{0.5\textwidth}
            \centering
            \begin{tikzonimage}[keepaspectratio, width=3 in]{"mte_duration_distance_to_along_optimal.png"}
            \node[fill=white, opacity=1.0, text opacity=1, anchor=south] at (0.5,-0.07) {$\hat{d}^S_T$};
            \node[fill=white, opacity=1.0, text opacity=1, anchor=south, rotate=90] at (0.05,0.5) {$\hat{d}^S_A$};
            \end{tikzonimage}
            \caption{Optimal Solutions}
            \label{fig: mte_duration_distance_to_along_optimal}
        \end{subfigure}
        \caption{
        $\hat{d}_T^\mathcal{S}$ and $\hat{d}_A^\mathcal{S}$ across all solution families are shown, with robust solutions categorized by $\Delta \tau$. The mean values for each solution category is shown with a marker, accompanied by their respective one standard-deviation covariance ellipsoids.
        }
        \label{fig: mte_duration_distance_to_along}
    \end{figure}
}
{
}

In summary, as shown in Fig. \ref{fig: mte_duration_distance_to_along}, the overall feasible solution set for the robust case exhibits higher $\hat{d}_T^\mathcal{S}$ and $\hat{d}_A^\mathcal{S}$ values compared to the non-robust case. 
However, if we only consider the optimal solutions, the mean $\hat{d}_T^\mathcal{S}$ for the robust solutions decreases, becoming more comparable to that of the non-robust solutions.
The shift in the mean between the feasible and optimal solutions suggests that, on average, the robust optimal solutions utilize the invariant manifolds as effectively as the non-robust optimal solutions. 
The covariance ellipsoids represent one standard deviation reveal a larger spread in the robust solutions, particularly in the $\hat{d}_A^\mathcal{S}$ direction, indicating greater variability in how these solutions flow along the invariant manifolds. 

\subsubsection{Dependence on $\tau_1$}
\label{subsubsection: results and discussion: quantitative analysis: dependence on tau}

In this section, we explore how the relationship of the robust solutions to the invariant manifolds changes as we vary $\tau_1$. 
We begin by grouping the solutions based on whether $\tau_1$ occurs on the forward shooting arc (first half of the shooting horizon) or the backward shooting arc (second half of the shooting horizon), and performing similar analysis to those discussed in \S \ref{subsubsection: results and discussion: quantitative analysis: dependence on delta tau}.
To distinguish between these two cases within the robust solution family, we will refer to the first case as `forward robust solutions', and the second as `backward robust solutions'.

\ifthenelse{\boolean{includefigures}}
{
    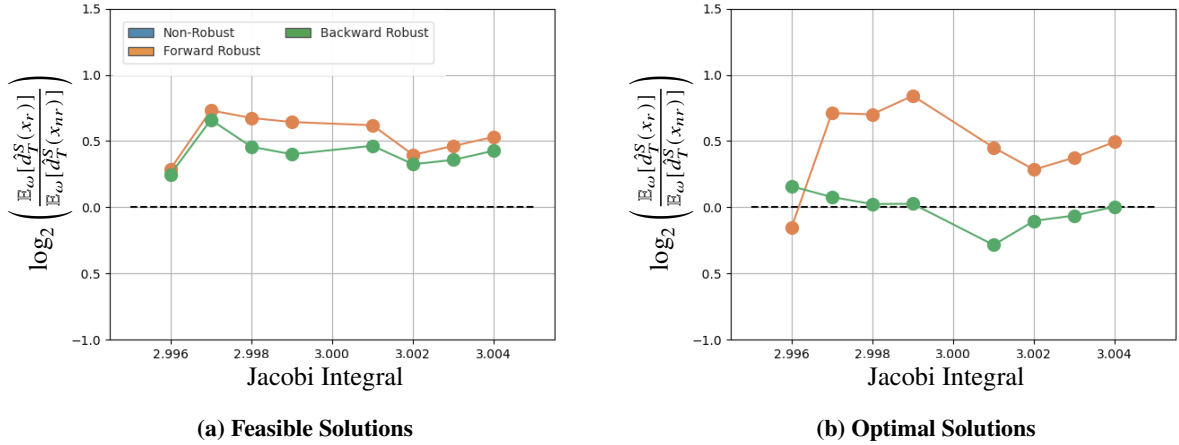
\begin{figure}[!htb]
        \centering
        \begin{subfigure}[b]{0.5\textwidth}
            \centering
            \begin{tikzonimage}[keepaspectratio, width=3 in]{"mte_start_distance_grouped_to_feasible.png"}
            \node[fill=white,opacity=0.85] at (0.42,0.80){\includegraphics[keepaspectratio, width=1.65 in]{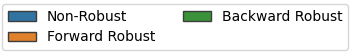}};
            \node[fill=white, opacity=1.0, text opacity=1, anchor=south] at (0.5,-0.03) {Jacobi Integral};
            \node[fill=white, opacity=1.0, text opacity=1, anchor=south, rotate=90] at (0.075,0.5) {$\text{log}_2 \left(\frac{\expectation{\omega}{\hat{d}^S_T (x_r)}}{\expectation{\omega}{\hat{d}^S_T (x_{nr})}}\right)$};
            \end{tikzonimage}
            \caption{Feasible Solutions}
            \label{fig: mte_start_distance_to_feasible}
        \end{subfigure}%
        \begin{subfigure}[b]{0.5\textwidth}
            \centering
            \begin{tikzonimage}[keepaspectratio, width=3 in]{"mte_start_distance_grouped_to_optimal.png"}
            \node[fill=white, opacity=1.0, text opacity=1, anchor=south] at (0.5,-0.03) {Jacobi Integral};
            \node[fill=white, opacity=1.0, text opacity=1, anchor=south, rotate=90] at (0.075,0.5) {$\text{log}_2 \left(\frac{\expectation{\omega}{\hat{d}^S_T (x_r)}}{\expectation{\omega}{\hat{d}^S_T (x_{nr})}}\right)$};
            \end{tikzonimage}
            \caption{Optimal Solutions}
            \label{fig: mte_start_distance_to_optimal}
        \end{subfigure}
        \caption{
        Fold change in $\hat{d}_T^\mathcal{S}$ for robust solutions with varying $\tau_1$ relative to non-robust solutions
        }
        \label{fig: mte_start_distance_to}
    \end{figure}
}
{
}

For both the forward and the backward categories, the robust feasible solutions exhibit higher $\hat{d}_T^\mathcal{S}$ values compared to the non-robust feasible solutions, with the backward robust solutions showing smaller distances than the forward robust solutions (Fig. \ref{fig: mte_start_distance_to_feasible}). 
However, when examining the optimal solutions, we observe that the forward robust solutions not only demonstrate a lower mean distance at the first energy level compared to the non-robust solutions but also a lower mean distance than the backward robust solutions (Fig. \ref{fig: mte_start_distance_to_optimal}).
Conversely, as the energy levels increase, the distance increases for the forward robust solutions, and decreases for the backward robust solutions. 
If $\tau_1$ occurs during the forward shooting arc, we expect robust solutions to adjust their control parameters to more effectively leverage the invariant manifolds at the initial energy levels; conversely, if $\tau_1$ occurs during the backward shooting arc, the robust solutions are likely to make similar adjustments for the later energy levels.
These expectations are consistent with the trends shown in Fig. \ref{fig: mte_start_distance_to}.
For both the forward and the backward robust solutions, it is important to note the greater reliance on the invariant manifolds for the robust solutions compared to the non-robust solutions during the initial and final energy levels respectively. 
We know that the absolute $\hat{d}_T^\mathcal{S}$ for both robust and non-robust solutions decreases when considering the optimal solution subset - therefore, the fact that robust solutions shadow the invariant manifolds even more closely at certain energy levels highlight the greater reliance on the invariant manifolds for the robust optimal solutions dependent on where $\tau_1$ occurs.

\ifthenelse{\boolean{includefigures}}
{
    \begin{figure}[!htb]
        \centering
        \begin{subfigure}[b]{0.5\textwidth}
            \centering
            \begin{tikzonimage}[keepaspectratio, width=3 in]{"rob_vs_nonrob_distance_to_optimal_unstable_start.png"}
            \node[fill=white,opacity=0.85] at (0.42,0.79){\includegraphics[keepaspectratio, width=1.65 in]{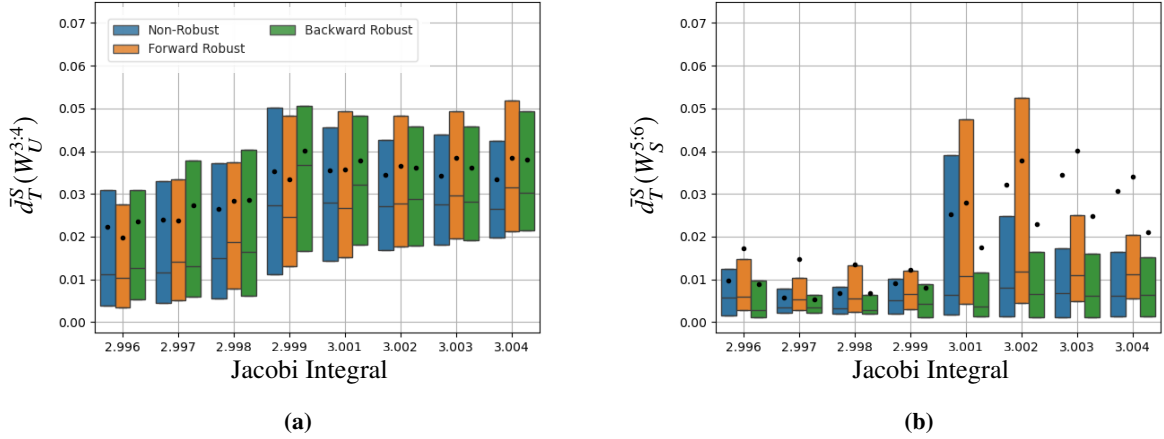}};
            \node[fill=white, opacity=1.0, text opacity=1, anchor=south] at (0.5,-0.03)  {Jacobi Integral};
            \node[fill=white, opacity=1.0, text opacity=1, anchor=south, rotate=90] at (0.05,0.5) {$\bar{d}_T^S(W_U^{3:4})$};
            \end{tikzonimage}
            \caption{}
            \label{fig: rob_vs_nonrob_distance_to_optimal_unstable_start}
        \end{subfigure}%
        \begin{subfigure}[b]{0.5\textwidth}
            \centering
            \begin{tikzonimage}[keepaspectratio, width=3 in]{"rob_vs_nonrob_distance_to_optimal_stable_start.png"}
            \node[fill=white, opacity=1.0, text opacity=1, anchor=south] at (0.5,-0.03)  {Jacobi Integral};
            \node[fill=white, opacity=1.0, text opacity=1, anchor=south, rotate=90] at (0.05,0.5) {$\bar{d}_T^S(W_S^{5:6})$};
            \end{tikzonimage}
            \caption{}
            \label{fig: rob_vs_nonrob_distance_to_optimal_stable_start}
        \end{subfigure}
        \caption{
        $\bar{d}_T^S(W_U^{3:4})$ and $\bar{d}_T^S(W_S^{5:6})$ for optimal non-robust solutions and robust solutions with varying $\tau_1$
        }
        \label{fig: rob_vs_nonrob_distance_to_optimal_start}
    \end{figure}
}
{
}

$\bar{d}_T^S(W_U^{3:4})$ increases across solution categories as we move through the energy levels (Fig. \ref{fig: rob_vs_nonrob_distance_to_optimal_unstable_start}). 
If we focus on the forward robust solutions, we observe that these solutions have a slightly smaller $\bar{d}_T^S(W_U^{3:4})$ compared to the non-robust solutions during the first energy level. 
This is expected, as robust solutions that experience an MTE early in the trajectory are likely to leverage $W_U^{3:4}$ more closely to compensate for the thruster outage.
A similar pattern is observed for $\bar{d}_T^S(W_S^{5:6})$ in robust backward solutions. 
These solutions are expected to align more closely with $W_S^{5:6}$, especially in the later energy levels, to mitigate the effects of a thruster outage in the latter half of the trajectory. 
Interestingly, these solutions maintain close alignment with $W_S^{5:6}$ across all energy levels, which is logical because compensating for a thruster outage in the latter half requires the trajectory to remain close to $W_S^{5:6}$ throughout, not just in the later stages (Fig. \ref{fig: rob_vs_nonrob_distance_to_optimal_stable_start}).

\ifthenelse{\boolean{includefigures}}
{
    \begin{figure}[!htb]
        \centering
        \begin{subfigure}[b]{0.5\textwidth}
            \centering
            \begin{tikzonimage}[keepaspectratio, width=3 in]{"rob_vs_nonrob_distance_along_optimal_unstable_start.png"}
            \node[fill=white,opacity=0.85] at (0.42,0.79){\includegraphics[keepaspectratio, width=1.65 in]{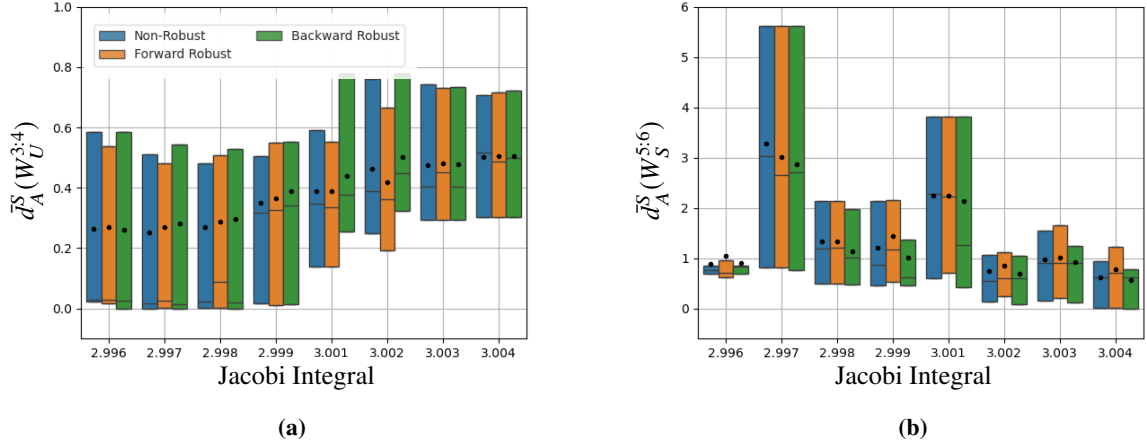}};
            \node[fill=white, opacity=1.0, text opacity=1, anchor=south] at (0.5,-0.03)  {Jacobi Integral};
            \node[fill=white, opacity=1.0, text opacity=1, anchor=south, rotate=90] at (0.075,0.5) {$\bar{d}_A^S(W_U^{3:4})$};
            \end{tikzonimage}
            \caption{}
            \label{fig: rob_vs_nonrob_distance_along_optimal_unstable_start}
        \end{subfigure}%
        \begin{subfigure}[b]{0.5\textwidth}
            \centering
            \begin{tikzonimage}[keepaspectratio, width=3 in]{"rob_vs_nonrob_distance_along_optimal_stable_start.png"}
            \node[fill=white, opacity=1.0, text opacity=1, anchor=south] at (0.5,-0.03)  {Jacobi Integral};
            \node[fill=white, opacity=1.0, text opacity=1, anchor=south, rotate=90] at (0.09,0.5) {$\bar{d}_A^S (W_S^{5:6})$};
            \end{tikzonimage}
            \caption{}
            \label{fig: rob_vs_nonrob_distance_along_optimal_stable_start}
        \end{subfigure}
        \caption{
        $\bar{d}_A^S(W_U^{3:4})$ and $\bar{d}_A^S (W_S^{5:6})$ for optimal non-robust solutions and robust solutions with varying $\tau_1$
        }
        \label{fig: rob_vs_nonrob_distance_along_optimal_start}
    \end{figure}
}
{
}

If we examine $\bar{d}_A^S(W_U^{3:4})$ and $\bar{d}_A^S(W_S^{5:6})$ for the robust solution categories, we notice some interesting trends in the solutions. 
We do not notice any discernible pattern for the robust solutions with respect to $\bar{d}_A^S(W_U^{3:4})$, but we notice that the robust backward solutions demonstrate a smaller $\bar{d}_A^S(W_S^{5:6})$ compared to robust forward solutions, as well as the non-robust solutions (Fig. \ref{fig: rob_vs_nonrob_distance_along_optimal_start}). 
This makes sense, as these backward arc solutions are expected to leverage the stable manifold $W_S^{5:6}$ more closely, especially in the latter half of the trajectory.
Interestingly, the robust forward solutions also maintain relatively small $\bar{d}_A^S(W_S^{5:6})$ compared to the non-robust solutions, especially toward the initial energy levels, suggesting that the robust forward solutions rely on $W_S^{5:6}$, similar to the robust backward solutions.

\ifthenelse{\boolean{includefigures}}
{
    \begin{figure}[!htb]
        \centering
        \begin{subfigure}[b]{0.5\textwidth}
            \centering
            \begin{tikzonimage}[keepaspectratio, width=3 in]{"mte_start_distance_to_along_feasible.png"}
            \node[fill=white,opacity=0.85] at (0.42,0.19){\includegraphics[keepaspectratio, width=1.65 in]{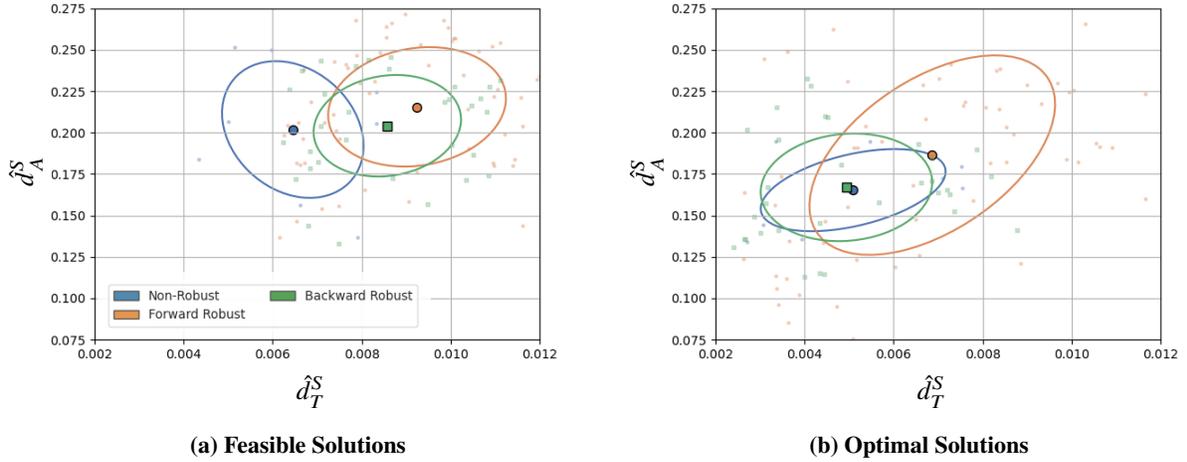}};
            \node[fill=white, opacity=1.0, text opacity=1, anchor=south] at (0.5,-0.07) {$\hat{d}^S_T$};
            \node[fill=white, opacity=1.0, text opacity=1, anchor=south, rotate=90] at (0.05,0.5) {$\hat{d}^S_A$};
            \end{tikzonimage}
            \caption{Feasible Solutions}
            \label{fig: mte_start_distance_to_along_feasible}
        \end{subfigure}%
        \begin{subfigure}[b]{0.5\textwidth}
            \centering
            \begin{tikzonimage}[keepaspectratio, width=3 in]{"mte_start_distance_to_along_optimal.png"}
            \node[fill=white, opacity=1.0, text opacity=1, anchor=south] at (0.5,-0.07) {$\hat{d}^S_T$};
            \node[fill=white, opacity=1.0, text opacity=1, anchor=south, rotate=90] at (0.05,0.5) {$\hat{d}^S_A$};
            \end{tikzonimage}
            \caption{Optimal Solutions}
            \label{fig: mte_start_distance_to_along_optimal}
        \end{subfigure}
        \caption{
        $\hat{d}_T^\mathcal{S}$ and $\hat{d}_A^\mathcal{S}$ across all solution families are shown, with robust solutions categorized by $\tau_1$.
        }
        \label{fig: mte_start_distance_to_along}
    \end{figure}
}
{
}

As shown in Fig. \ref{fig: mte_start_distance_to_along}, the average $\hat{d}_T^\mathcal{S}$ and $\hat{d}_A^\mathcal{S}$ is lower for the optimal solution subset compared to the overall feasible solution set (Fig. \ref{fig: mte_start_distance_to_along}).
In the backward case, the robust optimal solutions behave more similarly to the non-robust solutions.
Because we expect the dynamics to be more chaotic during the backward shooting arc, it becomes necessary for the robust backward solutions to shadow the invariant manifolds more closely.
Although $\hat{d}_T^\mathcal{S}$ and $\hat{d}_A^\mathcal{S}$ decreases with the optimal solutions in relation to the feasible solutions, they remain higher for the robust forward case compared to the non-robust solutions.
Because the forward shooting segments are farther away from Europa, it is not as crucial for the robust forward solutions to leverage the invariant manifolds as strongly as the robust backward solutions.

%% file: sections/conclusion.tex
\section{Conclusion}
\label{section: conclusion}

In this paper, we provide a comprehensive statistical analysis to compare robust solutions to non-robust solutions with respect to the underlying dynamical structures in a multibody gravitational environment.
We first place the missed thrust design problem in a general robust optimal control framework that can account for various forms of uncertainty.  
After this broad mathematical definition, we narrow our scope to the specific restricted missed thrust design problem and provide more detail on the specific mathematical formulation we use in this study. 
We compare a qualitative comparison between a robust optimal and non-robust optimal solution with respect to the invariant manifolds, and the quantitative behavior for a family of solutions using two distance metrics on a Poincar\'{e} section.
We illustrate the difference in the behavior of the feasible solutions and the optimal solutions for both non-robust and robust solutions.
We also present the difference in the behavior of the robust solutions depending on where the outage occurs, and its duration.
Our findings indicate that the robust optimal solutions align themselves to the pertinent invariant manifolds as closely as their non-robust counterparts, and in some cases demonstrate an even stronger alignment than the non-robust solutions.
Future efforts should investigate how the behavior will change as we lower the available control authority (e.g., lower thrust acceleration), and consider cases where the reference solution is coupled with multiple realization solutions. 
For both cases, we expect solutions to show stronger alignment with the invariant manifolds, and these trends becoming more prominent.
Having knowledge of the change in relation of a non-robust to a robust trajectory with respect to the pertinent dynamical structures should enable more efficient solution methods to robust trajectory problems, and in particular better enable the solution of these problems within a global optimization framework. 

%% file: sections/acknowledgement.tex
\section*{Acknowledgments}
The simulations presented in this article were performed on computational resources managed and supported by Princeton Research Computing, a consortium of groups including the Princeton Institute for Computational Science and Engineering (PICSciE) and the Office of Information Technology's High Performance Computing Center and Visualization Laboratory at Princeton University. The data that support this study are available from the authors upon reasonable request.